\documentclass{amsart}

\usepackage{amssymb,amsfonts, amscd,latexsym,amsmath,hhline,array,longtable}
\usepackage{pdfsync,color, comment,colortbl, mathrsfs,stmaryrd,cite,graphicx}
\usepackage[margin=10pt]{caption}

\usepackage{ifpdf}
\ifpdf \usepackage[colorlinks=true, citecolor=blue, linkcolor=blue, urlcolor=blue]{hyperref} \fi

\newcommand{\cal}{\mathcal}

\newtheorem{formula}{}[section]
\newtheorem{definition}[formula]{Definition}
\newtheorem{corollary}[formula]{Corollary}
\newtheorem{remark}[formula]{Remark}
\newtheorem{lemma}[formula]{Lemma}
\newtheorem{theorem}[formula]{Theorem}
\newtheorem{proposition}[formula]{Proposition}

\def\thrm{\begin{theorem}}
\def\thrml#1{\begin{theorem}\label{#1}}
\def\ethrm{\end{theorem}}
\def\rmrk{\begin{remark}}
\def\rmrkl#1{\begin{remark}\label{#1}}
\def\ermrk{\end{remark}}
\def\dfntn{\begin{definition}}
\def\dfntnl#1{\begin{definition}\label{#1}}
\def\edfntn{\end{definition}}
\def\nmrt{\begin{enumerate}}
\def\enmrt{\end{enumerate}}
\def\tm#1{\item[{\rm (#1)}]}
\def\qtnl#1{\begin{equation}\label{#1}}
\def\eqtn{\end{equation}}
\def\lmm{\begin{lemma}}
\def\lmml#1{\begin{lemma}\label{#1}}
\def\elmm{\end{lemma}}
\def\crllr{\begin{corollary}}
\def\crllrl#1{\begin{corollary}\label{#1}}
\def\ecrllr{\end{corollary}}
\def\css{\begin{cases}}
\def\ecss{\end{cases}}
\def\prf{\begin{proof}}
\def\eprf{\end{proof}}
\def\prpstn{\begin{proposition}}
\def\prpstnl#1{\begin{proposition}\label{#1}}
\def\eprpstn{\end{proposition}}

\def\cG{{\cal G}}

\def\cP{{\cal P}}

\def\cS{{\cal S}}

\def\cX{{\cal X}}
\def\cY{{\cal Y}}

\def\fc{{\frak c}}
\def\fC{{\frak C}}

\def\fF{{\frak F}}

\def\fX{{\frak X}}
\def\fY{{\frak Y}}

\def\sC{{\mathscr{C}}}

\DeclareMathOperator{\aut}{Aut}

\DeclareMathOperator{\diag}{Diag}

\DeclareMathOperator{\id}{id}
\DeclareMathOperator{\im}{im}

\DeclareMathOperator{\inv}{Inv}

\DeclareMathOperator{\mon}{Mon}
\DeclareMathOperator{\orb}{Orb}

\DeclareMathOperator{\pr}{pr}
\DeclareMathOperator{\res}{res}

\DeclareMathOperator{\rad}{Rad}

\DeclareMathOperator{\rk}{rk}

\DeclareMathOperator{\sym}{Sym}
\DeclareMathOperator{\WL}{WL}
\DeclareMathOperator{\WLD}{W_{1/2}}
\DeclareMathOperator{\Ws}{W}

\def\bone{{\bf 1}}
\def\bull{\hfill$\square$\medskip}
\def\grp#1{\langle {#1}\rangle}

\def\mltst#1{\{\hspace{-3pt}\{#1\}\hspace{-3pt}\}}
\def\qaq{\quad\text{and}\quad}

\def\phmaa#1{{\phantom{x}\hspace{-2mm}^{#1}}}
\def\simn#1{\stackrel{\scriptscriptstyle{_#1}}{\sim}}

\def\wh{\widehat}
\def\wt{\tilde}

\begin{document}

%\begin{comment}
\title{On the Weisfeiler algorithm of depth-$1$ stabilization}
\author{Gang Chen}
\address{School of Mathematics and Statistics, and Key Laboratory of Nonlinear Analysis $\&$ Applications (Ministry of Education), Central China Normal University, Wuhan 430079, China}
\email{chengangmath@mail.ccnu.edu.cn}
\author{Qing Ren}
\address{School of Mathematics and Statistics, and Key Laboratory of Nonlinear Analysis $\&$ Applications (Ministry of Education), Central China Normal University, Wuhan 430079, China}
\email{renqing@mails.ccnu.edu.cn}
\author{Ilia Ponomarenko}
\address{Steklov Institute of Mathematics at St. Petersburg, Russia; School of Mathematics and Statistics, and Key Laboratory of Nonlinear Analysis $\&$ Applications (Ministry of Education), Central China Normal University, Wuhan 430079, China}
\email{inp@pdmi.ras.ru}
\thanks{}
\date{}

\begin{abstract}
An origin of the multidimensional Weisfeiler-Leman algorithm goes back to a refinement procedure of deep stabilization,  introduced  by B.~Weisfeiler in a paper  included in the collective monograph ``On construction and identification of graphs'' (1976). This procedure is recursive and the recursion starts from an algorithm of depth-$1$ stabilization, which has never been discussed in the literature.

A goal of the present paper is to show  that a simplified algorithm of the depth-$1$ stabilization has the same power as the $3$-dimensional Weisfeiler-Leman algorithm. It is proved that the class of coherent configurations obtained at the output of this simplified algorithm coincides with the class introduced earlier by the third author. As an application we also prove that if there exist  at least two  nonisomorphic projective planes of order~$q$, then  the Weisfeiler-Leman dimension of the incidence graph of any projective plane of order~$q$ is at least $4$.
\end{abstract} 

\maketitle

\section{Introduction}\label{181023a}
In 1968, B.~Weisfeiler and A.~Leman proposed an algorithm $\WL$ for testing graph isomorphism \cite{WLe68}. The idea of the $\WL$ is as follows: analyzing the triples of vertices of a graph $\cG$ with vertex set~$\Omega$, construct a canonical partition $\WL(\cG)$ of the Cartesian square~$\Omega^2$, such that pairs from different classes lie in different orbits of the induced action of the group $\aut(\cG)$ on $\Omega^2$.  It should be noted that the goal (as Weisfeiler later wrote) was to construct an automorphic partition of~$\cG$, i.e., a partition of $\Omega$ into the orbits of the group $\aut(\cG)$. A year later, it became clear that the constructed algorithm did not achieve this goal~\cite{AdelWLF1969}. 

It can be seen from the 1968 paper that the authors understood the need to include individualization of vertices in the ``right'' algorithm. Probably, this was fully realized  in the paper by Weisfeiler, included by  him in the collective monograph~\cite[Section~O]{Wei1976}, where he introduced the concept of deep stabilization. Again, the goal was to find the automorphic partition of~$\cG$, but here the individualization of the vertices of~$\cG$ has already explicitly been used. The number of individualized vertices is just the depth of the algorithm $\Ws_m$ of \emph{depth-$m$ stabilization} proposed by Weisfeiler. The algorithm starts with $\WL$ and calls the algorithm~ $W_{m-1}$ ($m\ge 2$) recursively, followed by a depth-$1$ stabilization (the exact definitions can also be found in Section~\ref{010123h}). It was conjectured that if $m$ is at least $O(\log n)$, then the depth-$m$ stabilization determines the automorphic partition of
any graph of order~$n$. Furthermore, it was proved that depth-$1$ stabilization determines this partition in some nontrivial special cases~\cite[Section~P]{Wei1976}.

The algorithm now known as the multidimensional Weisfeiler-Leman algorithm was used by many
researchers in the late 1970s and early 1980s, see~\cite{Babai2015c}. The $m$-dimensional version $\WL_m$ of this algorithm (the \emph{$m$-dim~WL}) analyzes $(m+1)$-tuples of vertices, but the algorithm itself is non-recursive and does not explicitly use individualization of the vertices of input graph (note that $\WL_2=\WL$). This algorithm has been extensively studied in the context of the graph isomorphism problem and descriptive complexity theory (see~\cite{Grohe2017}) and culminated in its use in the Babai quasipolynomial algorithm~\cite{Babai2015c}. At the same time, it became clear that the $\WL_m$ cannot determine the automorphic partition of some graphs~$\cG$ even if $m=cn$ for an absolute constant~$c$ \cite{CaiFI1992} (see also~\cite{EvdP1999c}). 

One of the main goals of the present paper is to compare the power of the depth-$1$ stabilization $\Ws=\Ws_1$  and the $m$-dim WL. It should be noted that these algorithms being applied to a graph~$\cG$ have different outputs (for $m\ge 3$). The algorithm $\Ws$  constructs a partition $\Ws(\cG)$ of the set $\Omega^2$, that is a \emph{coherent configuration} in the sense of D.~Higman~\cite{Hig1975}  (or stationary graph in terms of B.~Weisfeiler and~A.~Leman), whereas the second one constructs the partition $\WL_m(\cG)$ of the Cartesian power~$\Omega^m$, which is an \emph{$m$-ary coherent configuration}\footnote{For $m=2$, this is an ordinary coherent configuration.} in terms of~L.~Babai~\cite{Babai2015c} (see also~\cite{AndresHelfgott2017}). A natural way to compare the outputs is to consider the projection $\pr_2 \WL_m(\cG)$ of the $m$-ary coherent configuration $\WL_m(\cG)$ to the first two coordinates, which is an ordinary coherent configuration. Our first result states that in this sense the power of $\Ws$ is between the powers of $\WL_3$ and $\WL_4$, more exactly, 
\qtnl{140323b}
\pr_2 \WL_3(\cG)\le \Ws(\cG)\le \pr_2 \WL_4(\cG) 
\eqtn
for every colored graph~$\cG$. The proof will be  given in Section~\ref{300323a}.

The exact definition of the depth-$1$ stabilization is quite cumbersome, see Section~\ref{010123h}. Therefore, we consider also a ``light'' version of the depth-$1$ stabilization, which was actually studied  by B.~Weisfeiler in more detail, and the output of which is also a coherent configuration (in general, smaller than that in the full version of the depth-$1$ stabilization), denoted below by $\WLD(\cG)$. Roughly speaking, $\WLD(\cG)$ can be thought as the output of the Weisfeiler-Leman algorithm applied to the graph~$\cG$ in which any pair of vertices is colored by a multiset of colors of this pair in each coherent configuration $\WL(\cG_\alpha)$, $\alpha\in\Omega$, where $\cG_\alpha$  is a colored graph obtained by individualization of~$\alpha$ in~$\cG$.

\thrml{050601a}
For every colored graph $\cG$, we have 
\qtnl{230323a}
\WL_3(\cG)=\WL_3(\WLD(\cG)) \qaq \pr_2 \WL_3(\cG)=\WLD(\cG).
\eqtn
\ethrm

Although our informal definition of $\WLD(\cG)$ is given for colored graphs, the concept of coloring is not necessary. Indeed, the coherent configurations which are stable with respect to~$\WLD$ are \emph{sesquiclosed} in the sense of~\cite{Ponomarenko2022a}, and, moreover, $\WLD(\cG)$ is just the intersection of all coherent configurations on~$\Omega$, containing the edge set of~$\cG$ as one of the relations (Theorem~\ref{041222b}). In fact, the two formulas in Theorem~\ref{050601a} express some relations  between $3$-ary and $2$-ary coherent configurations.

On the other hand, the concept of coloring is essential if one is  interested  in applying the depth-$m$ stabilization and the $m$-dim $\WL$ to test isomorphism directly, without computing an automorphic partition. In this case, it is said that (a) graphs~$\cG$ and $\cG'$ are \emph{$\WL_m$-equivalent} (or  not distinguished by the algorithm $\WL_m$) if the colorings of the partitions $\WL_m(\cG)$ and $\WL_m(\cG')$, constructed by the $\WL_m$ are the same, and (b) the~$\WL_m$ \emph{identifies} a graph $\cG$ if every graph $\WL_m$-equivalent to~$\cG$ is isomorphic to~$\cG$. In  a similar way, one can define the $\WLD$-equivalence and $\WLD$-identification.

\thrml{130123b}
Two  graphs  are $\WL_3$-equivalent if and only if they are $\WLD$-equiva\-lent. Moreover, the $\WL_3$ and $\WLD$  identify the same graphs.
\ethrm

The structure of the colorings used in the above definitions is quite complicated, for example, the coloring constructed by the $\WL_m$  contains full information about the intersection numbers of the corresponding $m$-ary coherent configuration. From  algebraic point of view, the colorings of two $m$-ary coherent configurations are equal if and only if the color preserving bijection between their color classes is an algebraic isomorphism between these coherent configurations  (for $m=2$,  an algebraic isomorphism is what  was called a weak equivalence in~\cite{Wei1976}). The proof of Theorem~\ref{130123b} is given just in terms of the algebraic isomorphisms. In this language, the algorithm $\WL$ identifies a graph $\cG$ if and only if the coherent configuration~$\WL(\cG)$ is separable, i.e., every algebraic isomorphism from it is induced by a (combinatorial) isomorphism~\cite[Theorem~2.5]{Fuhlbr2018a}. In particular, it gives a method to prove that a given graph has the \emph{Weisfeiler-Leman dimension} at most~$2$, where the Weisfeiler-Leman dimension of a graph $\cG$ is the smallest~$m$ for which the $\WL_m$ identifies~$\cG$.

In general, not every algebraic isomorphism   of a coherent configuration respects the property of being sesquiclosed. In this connection, it is quite natural to consider sesquiclosed algebraic isomorphisms, i.e., those that can be extended to algebraic isomorphisms of one-point extensions. The graphs  $\cG$ of the Weisfeiler-Leman dimension at most~$3$ can be characterized as those for which the coherent configuration~$\WLD(\cG)$ is \emph{sesquiseparable} in the sense that every sesquiclosed algebraic isomorphism from it is induced by an isomorphism. As an immediate consequence of Theorem~\ref{130123b} we obtain the following statement.
 
\crllrl{211222w}
The  Weisfeiler-Leman dimension of a graph $\cG$ is at most~$3$ if and only if the coherent configuration~$\WLD(\cG)$ is sesquiseparable.
\ecrllr

To illustrate the above concepts, we consider the incidence graphs of finite projective planes. On the first glance, the results obtained in~\cite{EvdP2010} make impression that the~$\WL_3$ cannot distinguish these graphs. Actually, this is true but a direct deduction does not give a nontrivial lower bound on the Weisfeiler-Leman dimension. On the other hand, the technique developed in the present paper enables us to do this.

\thrml{020123i}
The incidence  graphs of finite projective planes of the same order are $\WLD$-equivalent.
\ethrm

Theorems~\ref{020123i} and~\ref{050601a} show that any two incidence  graphs of finite projective planes of the same order~$q$ are $\WL_3$-equivalent. Thus the Weisfeiler-Leman dimension of any such graph is at least~$4$ if there are at least two nonisomorphic projective planes of  order~$q$ (which is the case for infinitely many prime powers~$q$).

\crllrl{230123i}
Let $q$ be a positive integer such that there are at least two nonisomorphic projective planes of order~$q$. Then the  Weisfeiler-Leman dimension of  the incidence graph of every finite projective plane of order~$q$  is at least~$4$.
\ecrllr

The paper is organized as follows. In Sections~\ref{050622g} and~\ref{120223a}, we give a necessary background on ordinary and $m$-ary coherent configurations,  and  the $m$-dimensional Weisfeiler-Leman algorithm. A general description of the depth-$m$ stabilization in the language of coherent configurations is given in Section~\ref{010123h}.  Special attention is paid here to the case $m=1$. Then in Section~\ref{041222}, we introduce  the $\WLD$-closure of a coherent configuration  and show that $\WLD$ is not stronger than~$\WL_3$.  Formula~\eqref{140323b} and Theorems~\ref{050601a} and~\ref{130123b} are proved in Section~\ref{300323a}. In Section~\ref{120223r}, we prove Theorem~\ref{020123i}.

\section{Coherent configurations}\label{050622g}

In order to make the paper as selfcontained as possible, we present  a necessary background from the theory of coherent configurations. Our notation is compatible with those in the monograph~\cite{CP2019}; proofs of the statements below can be found there.

\subsection{Notation}

Throughout the paper, $\Omega$ denotes a finite set. The $i$th coordinate of an $m$-tuple $x\in\Omega^m$ is denoted by $x_i$, $i=1,\ldots,m$.  The concatenation of $x$ and  $y\in\Omega^k$ is defined to be the $(m+k)$-tuple $x\cdot y=(x_1,\ldots,x_m,y_1,\ldots,y_k)$. %When $m=1$, we abbreviate $x_1\cdot y:=x\cdot y$. 

% Given a set $X$ of $m$-tuples, we put $X\cdot y=\{x\cdot y:\ x\in X\}$.

For $\Delta\subseteq \Omega$, the Cartesian product $\Delta\times\Delta$ and its diagonal are denoted by~$\bone_\Delta$ and $1_\Delta$, respectively; for $\delta\in\Delta$, we abbreviate  $1_\delta=1_{\{\delta\}}$.  For a relation $s\subseteq\bone_\Omega$, we set $s^*=\{(\alpha,\beta): (\beta,\alpha)\in s\}$, $\alpha s=\{\beta\in\Omega:\ (\alpha,\beta)\in s\}$ for all $\alpha\in\Omega$. %, and define $\grp{s}$ as the minimal (with respect to inclusion) equivalence relation on~$\Omega$, containing~$s$. 
For relations $r,s\subseteq\Omega\times\Omega$, we put $r\cdot s=\{(\alpha,\beta)\!:\ (\alpha,\gamma)\in r,\ (\gamma,\beta)\in s$ for some $\gamma\in\Omega\}$.
For any set $S$ of relations, we denote by $S^\cup$ the set of all unions of elements of~$S$. For any mapping $f:\Omega\to\Omega'$, we set $s^f=\{(\alpha^f, \beta^f):\ (\alpha, \beta)\in s\}$.

The partitions of the same set are partially ordered in a standard way: $S\le T$ if each class of $S$ is a union of some classes of~$T$, or equivalently, $S^\cup\subseteq T^\cup$. The join of $S$ and $T$ is defined as a uniquely determined partition $S\vee T$ such that $(S\vee T)^\cup=S^\cup\cap T^\cup$. For any set $\cS$ of partitions of the same set, that is closed with respect to the join, one can define a mapping $\fc$ from the set of all partitions of this set to~$\cS$, that takes   $S$  to the join of all $T\in\cS$ such that $T\ge S$. This mapping is a closure operator in the sense that for each partition $S$, we have $\fc(S)\ge S$, $\fc(\fc(S))=\fc(S)$, and $\fc(S)\ge \fc(T)$ if $S\ge T$ \cite{BurrisS2012}.

\begin{comment}
The set of classes of  an equivalence relation $e$ on~$\Omega$ is denoted by $\Omega/e$. For $\Delta\subseteq\Omega$, we set $\Delta/e=\Delta/e_\Delta$ where $e_\Delta=\bone_\Delta\cap e$. If the classes of $e_\Delta$ are singletons, $\Delta/e$ is identified with $\Delta$. Given a relation $s\subseteq \bone_\Omega$, we put
\qtnl{130622a}
s_{\Delta/e}=\{(\Gamma,\Gamma')\in \bone_{\Delta/E}: s^{}_{\Gamma,\Gamma'}\ne\varnothing\},
\eqtn
where $s^{}_{\Gamma,\Gamma'}=s\cap (\Gamma\times \Gamma')$. We also abbreviate $s^{}_\Gamma=s^{}_{\Gamma,\Gamma}$. 
Among all  equivalence relations $e$ on~$\Omega$, such that
\qtnl{kgu}
s=\bigcup_{(\Delta,\Delta')\in s_{\Delta/e}}\Delta\times \Delta',
\eqtn
there is  the largest (with respect to inclusion) one, which is denoted by  $\rad(s)$ and called the {\it radical} of~$s$. Obviously, $\rad(s)\subseteq\grp{s}$.
\end{comment}

\subsection{Rainbows}\label{130123a}
Let $S$ be a partition of $\Omega^2$. A pair $\cX=(\Omega,S)$ is called a \emph{rainbow} on $\Omega$ if 
\nmrt
\tm{C1}  $1_\Omega\in S^\cup$,
\tm{C2} $s^*\in S$ for all $s\in S$.
\enmrt
%The rainbow $\cX$ is said to be \emph{symmetric} if $s^*=s$ for all $s\in S$. 
The numbers $|\Omega|$ and $\rk(\cX)=|S|$ are called the {\it degree} and {\it rank} of~$\cX$, respectively. The elements of $S$ and of $S^\cup$ are called {\it basis relations} and \emph{relations} of~$\cX$. A unique basis relation containing the pair~$(\alpha,\beta)$ is denoted by~$r(\alpha,\beta)=r_\cX(\alpha,\beta)$.  The set of all relations of~$\cX$ is closed with respect to intersections and unions. Any $\Delta\subseteq\Omega$ for which $1_\Delta\in S$ is called a {\it fiber} of $\cX$; the set of all of them is denoted by $F=F(\cX)$. In view of condition~(C1), the set  $F$ forms a partition of~$\Omega$.  Any element of  $F^\cup$ is called a \emph{homogeneity set} of~$\cX$.

Let $\cX'=(\Omega', S')$ be a rainbow. A  {\it combinatorial isomorphism} or, briefly, \emph{isomorphism} from $\cX$ to $\cX'$ is defined to be a bijection $f: \Omega\rightarrow \Omega'$ such that  $s^f\in S'$ for all $s\in S$. In this case,  the rainbows $\cX$ and $\cX'$ are said to be {\it isomorphic}. The group of all isomorphisms of $\cX$ to itself contains a normal subgroup
$$
\aut(\cX)=\{f\in\sym(\Omega):\ s^f=s \text{ for all } s\in S\}
$$
called the {\it automorphism group} of $\cX$.

The partial order of the partitions of~$\Omega^2$ induces a partial order of the rainbows  on~$\Omega$. Namely, given two such rainbows~ $\cX$ and $\cY$, we set
$$
\cX\le\cY\ \Leftrightarrow\ S(\cX)^\cup\subseteq S(\cY)^\cup
$$
and say that $\cX$ is an \emph{fusion} of~$\cY$ and $\cY$ is an \emph{extension} of~$\cX$. The minimal and maximal elements with respect to this order are the {\it trivial} and {\it discrete} rainbows, respectively; in the former case, $S$ consists of~$1_\Omega$ and its complement (unless $\Omega$ is a singleton), and in the last case, $S$ consists of singletons. The \emph{intersection} of rainbows $\cX$ and $\cY$ is defined as  a uniquely determined rainbow $\cX\cap\cY$ such that $S(\cX\cap\cY)=S(\cX)\vee S(\cY)$, or equivalently,
$$
S(\cX\cap\cY)^\cup=S(\cX)^\cup\,\cap\, S(\cY)^\cup.
$$

A  {\it similarity}  from $\cX$ to $\cX'$ is defined to be a bijection $\varphi: S\rightarrow S'$ such that for all $s\in S$ we have $\varphi(s)\subseteq 1_{\Omega'}\ \Leftrightarrow\ s\subseteq 1_\Omega$, and  $\varphi(s^*)=\varphi(s)^*$.  In this case,  $\cX$ and~$\cX'$ are said to be {\it similar}. Extending $\varphi$ in a natural way to  bijections $S^\cup\to {S'}\phmaa{\cup}$ and $F^\cup\to {F'}\phmaa{\cup}$, denoted also by~$\varphi$, we have
 $$
\varphi(1_{\Omega^{}})=1_{\Omega'}\qaq \Omega^\varphi=\Omega'.
$$
Every isomorphism $f:\cX\to\cX'$ \emph{induces} the similarity $\varphi_f:\cX\to \cX',$ $s\mapsto s^f$. A similarity $\varphi$ is induced by an isomorphism $f$ if it holds that $\varphi=\varphi_f$. The identity similarity of~$\cX$ is denoted by~$\id_\cX$. A similarity~$\psi$ from a rainbow~$\cY\ge\cX$ \emph{extends} a  similarity~$\varphi$ from~$\cX$ if $\psi(s)=\varphi(s)$ for all $s\in S$.

\subsection{Coherent configurations}\label{190223j}
A rainbow $\cX$ is called a \emph{coherent configuration} if it satisfies the additional condition
\nmrt
\tm{C3} given $r,s,t\in S$, the \emph{intersection number} $c_{r,s}^t=|\alpha r\cap \beta s^{*}|$ does not depend on the choice of $(\alpha,\beta)\in t$. 
\enmrt
For any basis relation $t\in S$ such that $t\subseteq 1_\Omega$, the intersection number $c_{s,s^*}^t$  is either zero or is equal to $|\alpha s|$ for all $\alpha$ such that $\alpha s\ne\varnothing$;  in the latter case, the number $n_s=|\alpha s|$  is called the {\it valency} of~$s$. Every basis relation $s\in S$ is contained in the Cartesian product of two uniquely determined fibers which are called the \emph{left and right supports} of~$s$, respectively.  
 
 For  homogeneity sets $\Delta$ and $\Gamma$, we denote by~$S_{\Delta,\Gamma}$ the set of all  basis relations contained in $\Delta\times\Gamma$, and abbreviate $S_\Delta=S_{\Delta,\Delta}$. The rainbow $\cX_\Delta=(\Delta,S_\Delta)$ is a coherent configuration called the \emph{restriction} of~$\cX$ to~$\Delta$.   The coherent configuration~$\cX$ is said to be {\it homogeneous} or a {\it scheme} if $|F|=1$, or equivalently,  $1_\Omega\in S$.  

Let $G\le\sym(\Omega)$ be a permutation group. Denote by $S$ the set of all orbits $(\alpha,\beta)^G$ in the componentwise action of~$G$ on~$\Omega\times\Omega$, where $\alpha,\beta\in \Omega$. Then the pair 
$$
\inv(G)=\inv(G,\Omega)=(\Omega,S)
$$ 
is a coherent configuration. Any coherent configuration $\cX$ associated with some permutation group $G$ in this way is said to be {\it schurian}. In this case every  fiber $\Delta$ of~$\cX$  is an orbit of~$G$,  and the points of $\Delta$ can be identified with the right cosets of a point stabilizer $G_\delta$ with $\delta\in\Delta$,  so that the action of~$G$ on~$\Delta$ is induced by right multiplications. Note that a coherent configuration~$\cX$ is \emph{schurian} if and only if $\cX=\inv(\aut(\cX))$.  

Let $\cX$ and $\cX'$ be coherent configurations. A bijection $\varphi:S\rightarrow S'$, $r\mapsto r'$ is called an \emph{algebraic isomorphism} from $\cX$ to $\cX'$ if 
$$
c_{r', s'}^{t'}=c_{r^{},s^{}}^{t^{}}\quad\text{for all}\ r,s,t\in S.
$$
Every algebraic isomorphism of coherent configurations is also a similarity of the corresponding rainbows.  Furthermore, $\varphi$ preserves the degree, rank, and supports of basis relations, i.e., if $\Delta,\Gamma$ are fibers of $\cX$, then $r\in S_{\Delta,\Gamma}$ if and only if $r'\in S_{\Delta',\Gamma'}$.

Not every algebraic isomorphism is induced by an isomorphism. The coherent configuration $\cX$  is said to be \emph{separable} if every algebraic isomorphism from $\cX$  is induced by  an isomorphism.  Any trivial or discrete coherent configuration is separable. 

\subsection{Coherent closure}
The set of all coherent configurations on~$\Omega$ is closed with respect to intersection. The \emph{coherent closure} of a set $T\subseteq 2^{\Omega\times\Omega}$ is defined to be the intersection $\WL(T)$ of all coherent configurations~$\cX$ on $\Omega$, for which $T\subseteq S^\cup$. The notation reflects the fact that $\WL(T)$ can be computed by the algorithm $\WL$, see Section~\ref{181023a}.  

When $T$ is the union of the set $S(\cX)$ for some rainbow~$\cX$ on~$\Omega$ and a set $T'\subseteq 2^{\Omega\times\Omega}$, we put $\WL(\cX,T'):=\WL(T)$ and abbreviate $\WL(\cX):=\WL(\cX,\varnothing)$ and $\WL(\cX,s):=\WL(\cX,\{s\})$. In these notation, the mapping $\WL:\cX\mapsto\WL(\cX)$ becomes  a closure operator on the rainbows on~$\Omega$.

Let $\cX$ and $\cX'$ be rainbows and $\varphi:\cX\to\cX'$ a similarity. In general, there is no algebraic isomorphism $\psi:\WL(\cX)\to\WL(\cX')$ extending~$\varphi$. However, if such~$\psi$ does exist, then it is uniquely determined (and can be computed by the algorithm~$\WL$). Moreover, if $\varphi$ is induced by an isomorphism~$f$, then $f$ induces~$\psi$.

Let $m\ge 1$ and $y\in \Omega^m$. The  {\it $m$-point extension} of  a coherent configuration~$\cX$ on~$\Omega$ with respect to~$y$ is defined to be the coherent closure $\cX_y=\WL(\cX,T)$ with $T=\{1_{y_i}:\ 1\le i\le m\}$.  To simplify notation, we write $r_y$ instead of $r_{\cX_y}$. When the points are irrelevant, we use the term ``point extension'',  and if $m=1$, we abbreviate $\cX_{y_1}:=\cX_y$. 

Let $\varphi:\cX\to\cX'$ be an algebraic isomorphism, and let $x\in\Omega^m$, %$x'\in{\Omega'}^m$. \varphi_{xx'}
$x'\in\Omega'\phmaa{m}$. An  algebraic isomorphism $\psi:\cX^{}_{x^{}} \to\cX'_{x'}$ is  called an {\it $xx'$-extension } (or just a point extension) of $\varphi$ if $\psi$ extends $\varphi$  and
$$
\psi(1_{x^{}_i})=1_{x'_i},\quad i=1,\ldots,m.
$$
The $xx'$-extension is unique if it exists, and then we denote it by $\varphi_{xx'}$. Note that if the identity algebraic isomorphism $\varphi=\id_\cX$ has the extension $\varphi_{\wt xx}$ for some $m$-tuple~$\wt x$, then the  composition of $\varphi_{\wt xx}$ with $\varphi_{xx'}$ is obviously  the $\wt xx'$-extension of~$\varphi$. 

If $m=1$ and $y=(\alpha)$, $y'=(\alpha')$ are one-tuples,  then the algebraic isomorphism $\varphi_{\alpha\alpha'}:=\varphi_{yy'}$ is called  the $\alpha\alpha'$-extension of $\varphi$. Note that it  exists only if $\Delta^\varphi=\Delta'$, where $\Delta$ (respectively,~$\Delta'$) is the fiber of $\cX$ (respectively,~$\cX'$), containing the point~$\alpha$ (respectively,~$\alpha'$). 

When $\cX=\cX'$ and $\varphi=\id_\cX$, we introduce a binary relation $\sim$ on the concatenations $x\cdot y\in\Omega^{m+2}$, where $x\in\Omega^2$, so that
\qtnl{140323a1}
x\cdot y\sim x'\cdot y'\quad\Leftrightarrow\quad \varphi_{yy'}(r_{y^{}}(x))=r_{y'}(x');
\eqtn
in accordance with the above, here it is implicitly assumed that the  algebraic automorphism~$\varphi$ has the $yy'$-extension. 

\begin{comment}
An algebraic isomorphism $\varphi:\cX\to\cX'$ is said to be \emph{$\WL_m$}-closed if it has $(\alpha,\alpha')$-extension $\varphi_{\alpha,\alpha'}$ for any two compatible tuples  $\alpha\in\Omega^{m-1}$ and $\alpha'\in(\Omega')^{m-1}$. Every algebraic isomorphism is obviously $1$-closed.

A coherent configuration $\cX$ is said to be $\WL_m$-separable if for any coherent configuration $\cX'$, every $\WL_m$-closed  algebraic isomorphism $\varphi:\cX\to\cX'$ is induced by an isomorphism.  In particular, a  coherent  configuration is $\WL_1$-separarable if and only if it is separable. Denote by $s_{\scriptscriptstyle\WL}(\cX)$ the minimal integer $m$ for which $\cX$ is $\WL_m$-separable.

\prpstnl{050622d}
Let $G$ be a graph and $\cX=\WL(G)$. Then
$$
\dim_{\scriptscriptstyle\WL}(X)\le s_{\scriptscriptstyle\WL}(\cX)+2.
$$
\eprpstn
\end{comment}

\subsection{$m$-ary coherent configurations}\label{300323c}
Let  $m$ be a positive integer. For a tuple $x\in\Omega^m$, denote by $\rho(x)$ the equivalence relation on $M=\{1,\ldots,m\}$ such that $(i,j)\in\rho(x)$ if and only if  $x_i=x_j$.   The monoid of all maps from $M$ to itself is denoted by $\mon(M)$. Given $\sigma\in\mon(M)$, we set 
$$
x^\sigma=(x_{1^\sigma},\ldots,x_{m^\sigma}).
$$ 
%This notation is extended to  all tuples $y\in\Omega^{m-k}$ with $1\le k<m$ and all $\sigma$ with $\im(\sigma)=\{1,\ldots,m-k\}$ as follows: $y^\sigma=(y_{1^\sigma},\ldots,y_{k^\sigma})$.

Let $\fX$ be a partition of~$\Omega^m$ and $x\in\Omega^m$. The class of  $\fX$ containing $x$  is denoted by~$[x]$. Given $X_1,\ldots,X_m\in \fX$, we denote by $n(x;X_1,\ldots,X_m)$ the number of all $\alpha\in\Omega$ such that $x_{i\leftarrow \alpha} \in X_i$ for all~$i\in M$, where
\qtnl{030422q}
x_{i\leftarrow \alpha}=(x_1,\ldots,x_{i-1},\alpha,x_{i+1},\ldots,x_m).
\eqtn

A partition $\fX$ of  $\Omega^m$ is called an \emph{$m$-ary coherent configuration} on~$\Omega$  if the following conditions are satisfied for all $X\in \fX$:
\nmrt
\tm{C1'} $\rho(x)$ does not depend on $x\in X$,
\tm{C2'} $X^\sigma\in \fX$  for all $\sigma\in\mon(M)$,
\tm{C3'} for any $X_0,X_1,\ldots,X_m\in \fX$, the number $n_{X_1,\ldots,X_m}^{X_0}=n(x_0;X_1,\ldots,X_m)$ 
does not depend on $x_0\in X_0$.
\enmrt
The $2$-ary coherent configurations are the coherent configurations in the sense of Subsection~\ref{190223j}, and the definitions below are compatible with that introduced there. %In what follows the class of~$\fX$ containing $x\in\Omega^m$ is denoted by~$[x]$.

The $m$-ary coherent configurations on $\Omega$ are partially ordered in accordance with the partial order of the partitions of~$\Omega^m$. The intersection $\fX\cap\fY$ of two $m$-ary coherent configurations $\fX$ and $\fY$ on $\Omega$ is defined as the join of the corresponding partitions of~$\Omega^m$. It is not hard to verify that $\fX\cap\fY$ is also an $m$-ary coherent configuration. This enables us to define the $m$-ary coherent closure $\WL_m(\fX)$ for any partition~$\fX$ of~$\Omega^m$ as the intersection of all $m$-ary coherent configurations $\fY\ge\fX$. In particular, the mapping $\fX\mapsto\WL_m(\fX)$ is a closure operator, and  if $m=2$, then $\WL_m(\fX)=\WL(\fX)$.

Let $X\subseteq \Omega^m$, $k\in M$, and $y\in \Omega^{m-k}$. Denote by~$\pr_k X$ (respectively, $\res_yX$) the set of all $x\in\Omega^k$ such that $x\cdot x'\in X$ for some $x'\in\Omega^{m-k}$ (respectively, such that $x\cdot y\in X$). In accordance with \cite{AndresHelfgott2017} (see also \cite{Ponomarenko2022a}), if $\fX$ is an $m$-ary coherent configuration, then 
\qtnl{270123a}
\pr_k \fX=\{\pr_k X:\ X\in\fX\}\qaq \res_y\fX=\{\res_yX:\ X\in\fX,\ \res_yX\ne\varnothing\}
\eqtn
are  $k$-ary coherent configurations on~$\Omega$, called the \emph{$k$-projection} of~$\fX$,  and  the \emph{residue}  of~$\fX$ with respect~to~$y$, respectively. Obviously, 
\qtnl{180323a}
\res_y\fX\ge  \pr_k\fX.
\eqtn

For a set $X\in\fX$, the number of all $y\in X$ such that $\pr_k(x)=\pr_k(y)$ does not depend on~$x\in X$ \cite[Lemma~3.5]{Chen2023}, and is denoted by $n_k(X)$. Next, let  $y\in\Omega^{m-k}$, $t\in\{1,\ldots,m-k\}^k$, and  $y_{k,t}=(y_{t_1},\ldots,y_{t_k})$. Then it is easily seen that 
\qtnl{190323a}
\{y_{k,t}\}\in \res_y\fX.
\eqtn
 
Let $\fX'$ be an $m$-ary coherent configuration on~$\Omega'$. An \emph{isomorphism} from~$\fX$ to~$\fX'$ is  a bijection $f: \Omega^m\rightarrow {\Omega'}\phmaa{m}$ such that  $X^f\in \fX'$ and $(X^\sigma)^f=(X^f)^\sigma$ for all $X\in \fX$ and all $\sigma\in\mon(M)$. An \emph{algebraic isomorphism} from $\fX$ to $\fX'$ is a bijection $\varphi:\fX\to\fX'$ such that
\qtnl{110522a}
\varphi(X^\sigma)=\varphi(X)^\sigma\qaq n_{X_1,\ldots,X_m}^{X_0}=n_{\varphi(X_1),\ldots,\varphi(X_m)}^{\varphi(X_0)}
\eqtn
 for all $X,X_0,\ldots,X_m\in \fX$ and $\sigma\in\mon(M)$\footnote{The definition of algebraic isomorphisms given here is equivalent to Definition 7 in \cite{Chen2023}}. Again, for $m=2$, this definition is compatible with that given in Section~\ref{050622g}. For any $k\in M$, the algebraic isomorphism~$\varphi$ induces the  algebraic isomorphism between the $k$-projections,
 \qtnl{150223a}
 \varphi_k:\pr_k \fX\to\pr_k\fX',\ \pr_k X\mapsto \pr_k \varphi(X).
 \eqtn
 
 In a similar way, one can also define the $K$-projections $\pr_K \fX$  and $\varphi_K$ and the numbers $n_K(X)$ for any subset $K\subseteq M$. Any two tuples $y\in\Omega^k$ and $y'\in{\Omega'}\phmaa{k}$ are said to be \emph{$\varphi$-similar}, where $k=|K|$,  if the algebraic isomorphism $\varphi_K$ takes $[y]$ to~$[y']$.

\lmml{040622a}
Let $\varphi:\fX\to\fX'$ be an algebraic isomorphism of $m$-ary coherent configurations,  $k\in M$, and let $y$ and $y'$ be $\varphi$-similar $(m-k)$-tuples. Then  the mapping  
\qtnl{110423a}
\varphi_{yy'}:\res_{y^{}}\fX\to\res_{y'}\fX',\ \res_{y^{}} X\mapsto \res_{y'}\varphi(X)
\eqtn
is an algebraic isomorphism extending $\varphi_k$. Moreover, $\varphi_{yy'}$ takes $\{y^{}_{k,t}\}$ to~$\{y'_{k,t}\}$ for all $t\in\{1,\ldots,m-k\}^k$.
\elmm
\prf
The first part of the statement was proved in~ \cite[Lemma~3.6]{Ponomarenko2022a}. To prove the second  statement, let $x=y_{k,t}\cdot y$ and $X=[x]$. Then $\res_y X=\{y_{k,t}\}$, see formula~\eqref{190323a}. In accordance with~\cite[formula~(10)]{Ponomarenko2022a}, we have $\rho(X)=\rho(\varphi(X))$. Therefore,
$$
\varphi(X)=\{z'_{k,t}\cdot z':\ z'\in \pr_{K'} \varphi(X)\},
$$
where $K'=\{k+1,\ldots,m\}$. Since $y$ and $y'$ are $\varphi$-similar, we have $y'\in  \pr_{K'}\varphi(X)$. Therefore the set $\res_{y'}\varphi(X)$ contains $y'_{k,t}$ and hence coincides with the singleton~$\{y'_{k,t}\}$. Thus,
  $$
 \varphi_{yy'}( \{y_{k,t}\})=\varphi_{yy'}(\res_{y^{}}X)=\res_{y'}\varphi(X)=\{y'_{k,t}\},
$$ 
as required.
\eprf

\section{The $m$-dim Weifeiler-Leman algorithm}\label{120223a}

\subsection{The Weisfeiler-Leman algorithm}\label{040223u}
Let $\cX=(\Omega,S)$ be a rainbow and $\fC$ a linear ordered set. Any injective mapping $c:S\to \fC$ is called a \emph{coloring} of $\cX$. The element $c(\alpha,\beta):=c(r(\alpha,\beta))$ is  called the \emph{color} of the pair~$(\alpha,\beta)$ as well of the relation~$r(\alpha,\beta)$. A rainbow  equipped with a coloring is said to be  \emph{colored}.  An isomorphism $f:\cX\to\cX'$  of colored rainbows with colorings $c$ and~$c'$ is an ordinary isomorphism that  preserves the colorings, i.e., $c'(s^f)=c(s)$ for all~$s\in S$. 

Let $m\ge 2$. For a given colored rainbow~$\cX$, the $m$-dimensional Weisfei\-ler-Leman algorithm $\WL_m$ constructs  a certain coloring $c(m,\cX)$ of the set $\Omega^m$. At the first stage, an initial coloring $c_0=c_0(m,\cX)$ of $\Omega^m$ is constructed  from a coloring $c'$ defined by the following condition: given $x,y\in\Omega^m$,  we have $c'(x)=c'(y)$ if and only if $\rho(x)=\rho(y)$ and 
\qtnl{120622v}
c(x_i, x_j)=c(y_i, y_j)\quad \text{ for all } i,j\in M,
%c_0(x)=(c(x_1,x_1),,\ldots,c(x_1,x_m),\ldots, c(x_m,x_1),\ldots,c(x_m,x_m)).
\eqtn
where $c$ is the coloring of the rainbow~$\cX$.
Namely, the color $c_0(x)$ of  every $x$ is set to be the tuple $(c'(x^\sigma))_{\sigma\in\mon(M)}$. 

Starting from the second stage, the initial coloring is refined step by step. Namely, if $c_i$ is the coloring constructed at the $i$th step ($i\ge 0$), then the color of an $m$-tuple $x$ in the coloring $c_{i+1}$ is defined to be 
%the as follows: $c_{i+1}(x)=c_{i+1}(y)$ if and only if $c_i(x)=c_i(y)$ and
$$
c_{i+1}(x)=(c_i(x),\mltst{(c_i(x_{1\leftarrow \alpha}),\ldots,c_i(x_{m\leftarrow \alpha})):\ \alpha\in\Omega}),
$$
where $\mltst{\cdot}$ denotes a multiset. The algorithm stops when $|\im(c_i)|=|\im(c_{i+1})|$ and the final coloring $c(m,\cX)$ is set to be~$c_i$.

The color classes of the coloring $c(m,\cX)$ form  an $m$-ary coherent configuration denoted below by~$\WL_m(\cX)$. It is equal to the $m$-ary coherent closure $\WL_m(\fX_0)$  of the partition~$\fX_0$ consisting of the color classes of the coloring $c_0(m,\cX)$. It is important to note that this fact enables us to define the $m$-ary coherent configuration $\WL_m(\cX)$ also  for any uncolored rainbow~$\cX$. It was observed in~\cite{EvdP1999c} (see also~\cite{Dawar2020}) that the mapping 
$$
\pr_2 \WL_m:\cX\mapsto \pr_2 \WL_m(\cX)
$$  
is a closure operator taking the rainbows on~$\Omega$ to the coherent configurations on~$\Omega$, and, moreover, that for all $2\le k\le m$, we have
\qtnl{300323f}
\pr_2 \WL_k(\pr_2 \WL_m(\cX))=\pr_2 \WL_m(\cX).
\eqtn
\begin{comment}
\prf
Let $\fX=\WL_m(\cX)$. Formula~\eqref{120622v} implies that $\pr_2\fX\ge \cX$. Furthermore, it is not hard to verify that the classes of the initial coloring $c_0(m,\pr_2\fX)$ are contained in~$\fX^\cup$. Therefore,
$\WL_m(\pr_2\fX)=\fX$. Finally, if $\cX\ge \cY$, then $\fX_0\ge \fY_0$ implying by induction that $\fX\ge\fY$, where $\fY_0$ and  $\fY$ are the partitions into the color classes of $c_0(m,\cY)$ and $c(m,\cY)$, respectively. It follows that $\pr_2\fX\ge\pr_2\fY$.  
\eprf
\end{comment}
In what follows we also need the lemma below. 

\lmml{100622c1}
Let $\cX$ be a rainbow on $\Omega$, $m\ge 3$,  and $\fX=\WL_m(\cX)$. Then for any tuple $y\in \Omega^{m-2}$,
$$
\res_y\fX\ge (\pr_2 \fX)_y.
$$
\elmm
\prf
Formula~\eqref{180323a} for $k=2$ implies that  $\res_y\fX\ge\pr_2\fX\ge \cX$. Next, from formula~\eqref{190323a}, it follows that for all $i\in\{1,\ldots,m-k\}$, the singleton $\{1_{y_i}\}=\{y_{2,t_i}\}$ with $t_i=(i,i)$ is a class of the coherent configuration $\res_y\fX$. Thus, 
$$
\res_y\fX\ge \WL(\pr_2\fX, T_y)= (\pr_2 \fX)_y,
$$
where $T_y=\{1_{y_1},\ldots,1_{y_{m-2}}\}$. 
\eprf

The coherent configuration $\pr_2 \WL_m(\cX)$ can naturally be treated as a colored rainbow in which the color of $s\in S$ is defined to be the minimum of the colors $c(X)$, where $c=c(m,\cX)$  and $X$ is a class of $\WL_m(\cX)$ such that $\pr_2 X=s$. This colored rainbow is canonical in the sense that every color preserving isomorphism from $\cX$ to $\cX'$ is also a color preserving  isomorphism from $\pr_2 \WL_m(\cX)$  to $\pr_2 \WL_m(\cX')$ .

An uncolored rainbow $\cX$ is  said to be \emph{$\WL_m$-equivalent} to an uncolored rainbow~$\cX'$ with respect to a similarity  $\varphi:\cX\to\cX'$ if there exists an algebraic isomorphism  $\hat\varphi:\WL_m(\cX)\to\WL_m(\cX')$ such that the $2$-projection $\hat\varphi_2$ of $\hat\varphi$ extends~$\varphi$.  We say that the $m$-dim WL \emph{identifies} the rainbow~$\cX$ if  for every rainbow $\cX'$, every similarity  $\varphi:\cX\to\cX'$ such that  $\cX$ and $\cX'$ are $\WL_m$-equivalent, is induced by an isomorphism. The \emph{$\WL$-dimension} $\dim_{\scriptscriptstyle\WL}(\cX)$ of the rainbow~$\cX$ is defined to be the smallest positive integer $m$ such that the $m$-dim WL identifies~$\cX$.

\subsection{Rainbows and graphs} An undirected loopless graph $\cG$ can be treated as a colored rainbow. Indeed, let $\Omega$ and $E$ be the vertex and edge set of~$\cG$, respectively. Denote by $\cX(\cG)$ the colored rainbow on $\Omega$, that has exactly three basis relations, namely, $1_\Omega$, $E$, and  the edge set of the complement of~$\cG$, that are colored in the colors $0$, $1$, and~$2$, respectively. This can easily be extended to  vertex colored graphs~$\cG$ if  one replaces $1_\Omega$ by relations~$1_\Delta$, where $\Delta$ runs over the color classes of~$\cG$.

In accordance with above definitions, two vertex colored graphs $\cG$ and~$\cG'$ are isomorphic if and only if so are the colored rainbows $\cX=\cX(\cG)$ and $\cX'=\cX(\cG')$. Moreover, if~$\cG$ and~$\cG'$ are of the same order, then there is a unique color preserving mapping $\varphi:S(\cX)\to S(\cX')$, and $\varphi$  is a similarity from $\cX$ to~$\cX'$, called   a \emph{standard} one. 

In this terminology,  the $m$-dim WL does not distinguish the graphs $\cG$ and $\cG'$ in the sense of~\cite{Grohe2017} if and only if the (uncolored) rainbows $\cX$ and $\cX'$ are $\WL_m$-equivalent with respect to the standard similarity~$\varphi$. Indeed, the ``only if'' part is an immediate consequence of the fact that the mapping 
\qtnl{190223w}
\hat\varphi:\WL_m(\cX)\to\WL_m(\cX'),\ X\mapsto {c'}\phmaa{-1}(c(X))
\eqtn
is an algebraic isomorphism such that the $2$-projection $\hat\varphi_2$ of $\wh\varphi$ extends the algebraic isomorphism~$\varphi$, where $c=c(m,\cX)$ and $c'=c(m,\cX')$. Conversely, assume that the rainbows~$\cX$ and~$\cX'$ are $\WL_m$-equivalent with respect to~$\varphi$. Then there is an algebraic isomorphism $\hat\varphi:\WL_m(\cX)\to\WL_m(\cX')$ such that $\hat\varphi_2$ extends~$\varphi$. The last condition implies that $\hat\varphi(\fX^{}_0)=\fX'_0$ and the equality
$$
c^{}_i(X)=c'_i(\hat\varphi(X))
$$
holds for $i=0$ and all $X\in\fX_0$. By induction on $i$, this equality holds for all~$i\ge 0$. Therefore the $m$-dim WL does not distinguish the graphs $\cG$ and $\cG'$.

Let $\cG$ be a vertex colored graph. It is not difficult to verify that the $m$-dim WL \emph{identifies}~$\cG$ in the sense of~\cite{Grohe2017} if  and only if  the $m$-dim WL identifies the uncolored rainbow $\cX(\cG)$, and the $\WL$-dimension of $\cG$ equals $\dim_{\scriptscriptstyle\WL}(\cX)$.

\subsection{The Cai-F\"urer-Immerman theorem}\label{060223e}
Let us recall (very briefly) the definition of the $\sC_m$-pebbling game~\cite{CaiFI1992}. We describe it as a pebbling game $\sC_m(\varphi)$ for two rainbows $\cX$ and $\cX'$ and a similarity $\varphi:\cX\to \cX'$. 

There are two players, called Spoiler and Duplicator, and $m$ pairwise distinct pebbles, each given in duplicate.   The game consists of rounds and each round consists of two parts. At the first part, Spoiler chooses a set~$A'$ of points in~$\cX$ or in~$\cX'$. Duplicator responds with a set $A$ in the other rainbow, such that $|A|=|A'|$ (if this is impossible, then Duplicator loses). At the second part, Spoiler places one of the pebbles\footnote{It is allowed to move previously placed pebbles to other 	vertices and place more than one pebble on the same vertex.}  on a point in~$A$. Duplicator responds  by placing the copy of the pebble on some point of~$A'$. 
 
The configuration  after a round is determined by a bijection $f:\Delta\to\Delta'$, where $\Delta\subseteq\Omega$ (respectively, $\Delta'\subseteq\Omega' $) is the set of points  in $\cX$ (respectively, in~$\cX'$), covered by pebbles, and any two points $\alpha\in\Delta$ and $f(\alpha)\in\Delta'$ are covered by the copies of the same pebble. Duplicator wins the round if 
$$
(s_{\Delta^{}})^f=\varphi(s)_{\Delta'}
$$
for all $s\in S$ such that $s_\Delta\ne\varnothing$; in the other words, the bijection~$f$ induces the isomorphism between the rainbows $\cX^{}_{\Delta^{}}$ and $\cX'_{\Delta'}$, that respects the similarity~$\varphi$. Spoiler wins if Duplicator  does not win.

The game starts from an initial configuration (which  is considered as the configuration after zero rounds), i.e.,  a pair $(x,x')\in\Omega^k\times {\Omega'}\phmaa{k}$  for which $k\le m$ and the points~$x^{}_i$ and~$x'_i$ are covered with copies of the same pebble, $i=1,\ldots,k$. We say that Duplicator  has a winning strategy  for the game $\sC_m(\varphi)$  on $\cX$ and $\cX'$ with  initial configuration $(x,x')$ if, regardless of Spoiler's actions, Duplicator wins after any number of rounds (for appropriate selections of the sets~$A'$ and  points in~$A$).
\medskip

{\bf Notation.} {\it  $T^m_\varphi(\cX,\cX')$ is the set of all pairs $(x,x')\in\Omega^m\times {\Omega'}\phmaa{m}$ such that Duplicator has a winning strategy in the pebble game  $\sC_{m+1}(\varphi)$ on $\cX$ and $\cX'$ with initial configuration $(x,x')$.
}\smallskip

One of the main results proved in \cite{CaiFI1992} is the equivalence~$(1)\Leftrightarrow(3)$ in Theorem~5.2 there; it gives a characterization of the set  $T^m_\varphi(\cX,\cX')$ in terms of the final colors of the $m$-ary coherent configurations $\fX=\WL_m(\cX)$ and $\fX'=\WL_m(\cX')$. In our notation, it can be stated as follows:
\qtnl{281023a}
(x,x')\in T^m_\varphi(\cX,\cX')\quad\Leftrightarrow\quad c(m,\cX)(x)=c(m,\cX')(x')
\eqtn
for all $(x,x')\in\Omega^m\times {\Omega'}\phmaa{m}$. The lemma below is an easy consequence of this relation.

\lmml{131222a}
Let $\cX$ and $\cX'$ be rainbows and $\hat\varphi:\WL_m(\cX)\to\WL_m(\cX')$ a bijection such that $\varphi:=\pr_2(\hat\varphi)$ takes~$S$ to~$S'$.  Then
\nmrt
 \tm{1} if $\hat\varphi$ is an algebraic isomorphism, then $\hat\varphi([x])=[x']$ for all $x\in\Omega^m$ and all $x'\in{\Omega'}\phmaa{m}$ such that $(x,x')\in T^m_\varphi(\cX,\cX')$,
 \tm{2} if $[x]\times\hat\varphi([x])\subseteq T^m_\varphi(\cX,\cX')$ for all $x\in\Omega^m$, then $\hat\varphi$ is an algebraic isomorphism.
 \enmrt
\elmm
\prf
Assume that $\hat\varphi$ is an algebraic isomorphism. Then it is uniquely determined from the condition on its $2$-projection, and hence coincides with algebraic isomorphism~\eqref{190223w}. Thus the first statement follows from formula~\eqref{281023a}. To  prove the second statement, it suffices to note that the assumption and formula~\eqref{281023a} imply that the final colors of the classes $[x]$ and $\hat\varphi([x])$ are the same for all $x\in\Omega^m$.  Consequently, the mapping $\hat\varphi$ coincides with algebraic isomorphism~\eqref{190223w}.
\eprf

\section{Deep stabilization in the sense of Weisfeiler}\label{010123h}

\subsection{A general scheme}In this subsection, we discuss  the concept of the deep stabilization proposed by B.~Weisfeiler in~\cite[Section~O]{Wei1976}. In accordance with the remark made by him at page~VII, all definitions, constructions and algorithms of that section should be attributed to Weisfeiler himself. It should be noted that the terms ``graph'' and ``stationary graph'' used there  correspond  in our terminology to the structures  satisfying condition~(C1) and coherent configurations, respectively. In both cases the set $S$ of basis relations is assumed to be linear ordered and the corresponding order can naturally be interpreted as a coloring of~$\Omega^2$. 

One of primary goals for Weisfeiler was to find an efficicent algorithm that given a coherent configuration~$\cX$ constructs either a certificate for the equality $\cX=\inv(\aut(\cX))$ or a canonical coherent configuration $\cX'>\cX$ such that 
$$
F(\cX')=\orb(\aut(\cX));
$$ 
if the  algorithm was polynomial-time, then one would obviously solve the graph isomorphism problem. To this end, he proposed an algorithm $\Ws_m$ of deep stabilization that, given an integer $m\ge 1$, called the depth, recursively constructs a  canonical coherent configuration~$\Ws_m(\cX)\ge\cX$ the automorphism group of which coincides with  $\aut(\cX)$. Though it was not written explicitly, it is not hard to see that
$$
\cX\le \Ws_1(\cX)\le\ldots\le\Ws_m(\cX)=\Ws_{m+1}(\cX)=\ldots=\inv(\aut(\cX))
$$
for some $m\le n$.  It follows that the  stabilization of sufficiently large depth is sufficient to identify a given graph up to isomorphism.

The exact definition\footnote{It appears in Subsection~6.1 of~\cite[Chapter~O]{Wei1976}, the two other possible definitions (in Subsections~6.2 and 6.3) are not given in detail.} of the \emph{depth-$m$ stabilization} $\Ws_m$  is recursive and requires an auxiliary algorithm $\sigma$ that given a coherent configuration $\cX$ on~$\Omega$, and a family 
$$
%\fF=\{\cY^{\alpha}\}_{\alpha\in\Omega}
\fF=\{\cY^{\alpha}:\ \alpha\in\Omega\}
$$ 
of the coherent configurations~$\cY^{\alpha}$ on~$\Omega$, constructs a new coherent configuration $\sigma(\cX,\fF)\ge\cX$. The coherent configurations~$\cY^\alpha$ should be ``invariant'' in some natural sense to provide the equality  $\aut(\sigma(\cX,\fF))=\aut(\cX)$. In the case of the depth-$m$ stabilization, one takes $\cY^\alpha=\cX_\alpha$ at the first step of the recursion and defines  the coherent configuration $\Ws_m(\cX)$ as follows:
$$
\Ws_m(\cX)=\css
\sigma(\cX,\{\cX_\alpha:\ \alpha\in\Omega\})   &\text{if $m=1$,}\\
\sigma(\cX,\{\Ws_{m-1}(\cX_\alpha):\ \alpha\in\Omega\})   &\text{if $m>1$,}\\
\ecss
$$
where $\cX_\alpha$ is the $\alpha$-extension of~$\cX$ with canonically ordered basis relations. Thus the algorithm $\Ws=\Ws_1$ of  the \emph{depth-$1$ stabilization} consists in refining the coherent configuration $\cX$ with the help of the algorithm~$\sigma$ applied for the family~$\fF$ of the coherent configurations $\cY^{\alpha}=\cX_\alpha$.

\subsection{The depth-$1$ stabilization}
We begin with a notation. Given two  pairs  $x,y\in\Omega^2$, denote by $n_y(x)$ the number of all pairs $y'$ for which $x\cdot y\sim x\cdot y'$, see formula~\eqref{140323a1}. Clearly,  $n_y(x)\ge 1$. When $y=(\alpha,\alpha)$, we abbreviate $n_\alpha(x):=n_y(x)$.

Let $I=\{1,\ldots,4\}$ and $i\in I$. Based on the above notation, we define a binary relation $\simn{i}$ on the points (if $i=1$) or on the pairs of points  (if $i\ne 1$) of the coherent configuration~$\cX$:\medskip

\noindent(1)  $\alpha\simn{1}\alpha'\ \Leftrightarrow\ \id_\cX$ has the $\alpha\alpha'$-extension, \\
\noindent(2) $x\simn{2} x'\ \Leftrightarrow\ r(x)=r(x')$ and for each $\alpha$ there is~$\alpha'$ such that $n_{\alpha^{}}(x)=n_{\alpha'}(x')$,\\
\noindent(3) $x\simn{3}x'\ \Leftrightarrow\ \id_\cX$ has the $xx'$-extension, \\
\noindent(4) $x\simn{4} x'\ \Leftrightarrow\ x\simn{3} x'$ and for each~$y$ there is $y'$ such that $n_{y^{}}(x)=n_{y'}(x')$.\\

It should be noted that the equivalence relations $\simn{3}$ and $\simn{4}$ are finer than, respectively, $\simn{1}$ and $\simn{2}$ (take $x=(\alpha,\alpha)$ and $x'=(\alpha',\alpha')$).  The reason we consider four equivalence relations instead of two is to bring our interpretation of the Weisfeiler algorithm closer to the original algorithm.

It is easily seen that (a) $\simn{1}$ is an equivalence relation on $\Omega$, each class of which is contained in some fiber of~$\cX$, and (b) if $i\ne 1$, then $\simn{i}$ is an equivalence relation on~$\Omega^2$, each class of which is contained in some basis relation of~$\cX$. The set of classes of each equivalence relation $\simn{i}$ can naturally be treated as a (canonically) linear ordered family of binary relations (for $i=1$, one can take the corresponding reflexive relations). Denote this family by~$T_i$ and put
\qtnl{110323a}
\sigma_i(\cX)=\WL(\cX,T_i),\qquad i\in I.
\eqtn
It is not hard to see (and this was also done  in~\cite[Section~O]{Wei1976}) that  $\sigma_i(\cX)\ge\cX$ and $\aut(\cX)=\aut(\sigma_i(\cX))$ for each $i$. 

Now we are ready to  describe  the depth-$1$ stabilization. Namely, the algorithm~$\Ws$ consists in successive stabilization of the input coherent configuration~$\cX$ with the help of auxiliary  subroutines~$\sigma_i$, $i\in I$, computing the coherent closure~\eqref{110323a} by using ordinary Weisfeiler-Leman algorithm. More precisely, the coherent configuration $\Ws(\cX)$ is the output of the following procedure
\qtnl{010123k}
\text{{\bf while }} \cX<\sigma_i(\cX) \text{ for some } i\in I \text{ {\bf do} } \cX:=\sigma_i(\cX)\text{  {\bf od}}
\eqtn

It should be noted that our description differs somewhat from Weisfeiler's, since he represented  coherent configurations as matrices whose elements are (pairwise non-commuting) variables; in particular, the relation $\simn{i}$ was not calculated explicitly, but was presented  implicitly, in the form of equalities of the coefficients of some polynomials in these variables.

Some properties of the coherent configurations closed with respect to the depth-$1$ stabilization are proved in \cite[Chapters~O and~P]{Wei1976}. It was emphasized by Weisfeiler that all of them are already fulfilled for those $\cX$ for which $\sigma_1(\cX)=\sigma_2(\cX)=\cX$. A characterization of the coherent configurations satisfying the latter condition follows from Lemma~\ref{010123u1} and Theorem~\ref{130123b} to be proved in the subsequent sections.

Summarizing the above, we note a difference between the Weisfeiler deep stabilization algorithm and the modern multidimensional Weisfeiler-Leman algorithm. Namely, the first one does not construct an explicit coloring of the Cartesian power of the underlying set (although such a coloring can easily be defined, see  Section~\ref{120223e}). The increase in depth is achieved by using a recursion followed by the application of the  ordinary ($2$-dimensional) Weisfeiler-Leman algorithm. %Taking our results into account, one can say that the modern multidimensional Weisfeiler-Leman algorithm is a non-recursive version of the Weisfeiler deep stabilization.

\section{Sesquiclosed coherent configurations and sesquiclosure}\label{041222}

In accordance with \cite{Ponomarenko2022a}, a coherent configuration $\cX=(\Omega,S)$ is said to be \emph{sesquiclosed} if the following two conditions are satisfied for all $\alpha,\alpha'\in\Omega$:
\nmrt
\tm{S1} $F(\cX_\alpha)=\{\alpha s:\ s\in S,\ \alpha s\ne\varnothing\}$,
\tm{S2} if $\alpha$ and~$\alpha'$ belong to the same fiber of~$\cX$, then  $\alpha\simn{1}\alpha'$.
\enmrt
The conditions (S1) and (S2) are independent. Indeed, if  $\cG$ and $\cG'$ are the Payley and Piesert graphs (see~\cite{Peisert2001}) of the same sufficiently large order, then the coherent configuration of their disjoint union $\cG\cup\cG'$ is homogeneous. It satisfies condition~(S1), because both $\cG$ and~$\cG'$ are rank~$3$ graphs, and it does not satisfy condition~(S2), because any one-point extension of~$\WL(\cG)$ is separable (see the proof of \cite[Theorem~1.1]{Ponomarenko2020a}) and hence  the identity algebraic isomorphism $\id_\cX$ has no the $\alpha\alpha'$-extension, where $\alpha$ is a vertex of~$\cG$ and $\alpha'$ is a vertex of~$\cG'$.

Now, let $\cX$ be the  scheme of the Shrikhande graph, see e.g. \cite[Example 2.6.17]{CP2019}. This scheme has transitive automorphism group and hence satisfies condition~(S2). On the other hand, $\cX$ is a scheme of rank~$3$, and  a straightforward computation shows that any one-point extension of~$\cX$ has exactly $4$ fibers. Thus, $\cX$ cannot satisfy condition~(S1). 

\begin{comment}
The following statement shows that in any a sesquiclosed coherent configuration, the classes of the equivalence relation $\simn{2}$ are basis relations.

\lmml{110123a}
Let $\cX$ be  a sesquiclosed coherent configuration.  Then for every $\alpha$ and $s\in S(\cX_\alpha)$ the union of all $s^\varphi$ taken over all possible $(\alpha,\beta)$-extensions of $\id_\cX$, is a basis relation of~$\cX$.
\elmm
\end{comment}

\thrml{041222b}
The intersection of any two sesquiclosed coherent configurations is sesqui\-closed.
\ethrm
\prf
Let $\cY$ be the intersection of sesquiclosed coherent configurations $\cX$ and~$\cX'$. To prove that $\cY$ satisfies condition~(S2), let $\Delta\in F(\cY)$. Note that $\Delta$ is a union of some classes  of the equivalence relation $\simn{1}$ defined for the coherent configuration~$\cY$. Take any such class $\Delta'$. Then condition~(S2) is satisfied for all $\alpha,\alpha'\in\Delta'$, and we are done if $\Delta'=\Delta$. Suppose that $\Delta'\subsetneq \Delta$.  Then $\Delta'$ is not a homogeneity set of either~$\cX$ or~$\cX'$. By symmetry, we may assume that the former holds. Since $\Delta\in F(\cY)\subseteq F(\cX)^\cup$, there are a fiber $\Gamma\in F(\cX_\Delta)$ and points
$$
\alpha\in\Gamma\setminus\Delta' \qaq \alpha' \in \Gamma\cap\Delta'. 
$$
In view of condition~(S2),  the algebraic automorphism~$\id_\cX$ has the $\alpha\alpha' $-extension. Its restriction  to~$\cY_\alpha\le\cX_\alpha$ yields the $\alpha\alpha' $-extension of the algebraic automorphism~$\id_\cY$. It follows that  $\alpha\simn{1}\alpha' $ and hence $\alpha\in\Delta'$, a contradiction.  

Let $\alpha\in\Omega$ and $s\in S(\cY)$ be such that $\alpha s\ne\varnothing$. Then there are $\Delta,\Gamma\in F(\cY)$ such that  $\alpha\in\Delta$ and $s\subseteq \Delta\times \Gamma$. Since $\alpha s$ is a homogeneity set of the coherent configuration~$\cY_\alpha$, there is $\Gamma_\alpha\in F(\cY_\alpha)$ such that
\qtnl{051222s}
\Gamma_\alpha\subseteq \alpha s\subseteq \Gamma.
\eqtn

By the above paragraph, $\cY$ satisfies condition~(S2). Hence the algebraic automorphism~$\id_\cY$ has the $\alpha\beta$-extension $\varphi_\beta$ for each $\beta\in\Delta$ (this extension is identity if  $\beta=\alpha$). Then $\Gamma_\beta:=(\Gamma_\alpha)^{\varphi_\beta}$ is a fiber of~$\cY_\beta$.  In view of~\eqref{051222s}, we have 
\qtnl{061222a}
\Gamma_\beta=(\Gamma_\alpha)^{\varphi_\beta}\subseteq (\alpha s)^{\varphi_\beta}=
\beta s^{\varphi_\beta}=\beta s\subseteq \Gamma.
\eqtn

We claim that the union $r$ of all relations $\{\beta\}\times \Gamma_\beta$ with $\beta\in\Delta$, belongs to the set  $S^\cup\,\cap\, {S'}\phmaa{\cup}=S(\cY)^\cup$. Let us verify that $r\in S^\cup$; similarly, one can verify that $r\in {S'}\phmaa{\cup}$.  Since $\cX_\alpha\ge\cY_\alpha$,  there is~$r_\alpha\in S^\cup$ such that $\Gamma_\alpha=\alpha r_\alpha$. By condition~(S2),  the algebraic automorphism $\id_\cX$ has the $\alpha\beta$-extension $\varphi_{\alpha\beta}$ for each $\beta\in\Delta$. Since it extends $\varphi_\beta$, we have
$$
\Gamma_\beta=(\Gamma_\alpha)^{\varphi_\beta}=(\Gamma_\alpha)^{\varphi_{\alpha\beta}}=(\alpha r_\alpha)^{\varphi_{\alpha\beta}}=\beta r_\alpha.
$$
Thus $r$ is equal to the union of all relations $\{\beta\}\times \beta r_\alpha$, $\beta\in\Delta$, which is in its turn is equal to~$1_\Delta\cdot r_\alpha\in S^\cup$.

From the claim, it follows  that $r$ is a relation of~$\cY$. Moreover, $r\subseteq\Delta\times\Gamma$ by formula~\eqref{061222a}. Since $\alpha s\supseteq \alpha r$, we obtain $s=r$. Thus, 
$$
\alpha s=\alpha r=\Gamma_\alpha\in F(\cY_\alpha),
$$
as required.
\eprf

Let  $\cX$ be a rainbow on~$\Omega$. At least one coherent configuration larger than or equal to~$\cX$ is sesquiclosed, e.g., the discrete one. Denote by $\WLD(\cX)$ the intersection of all sesquiclosed coherent configurations larger than or equal to~$\cX$.\footnote{It will be clear from Lemma~\ref{010123u1} that the coherent configuration $\WLD(\cX)$ is obtained after applying ``half'' of the depth-$1$ stabilization algorithm~$\Ws$.} By Theorem~\ref{041222b}, the coherent configuration~$\WLD(\cX)$ is sesquiclosed and we call it the \emph{sesquiclosure} of~$\cX$. The mapping $\cX\mapsto \WLD(\cX)$ defines a closure operator on the set of all rainbows on the same set. The following lemma shows that the sesquiclosure is not more powerful than the closure defined by the $3$-dim~WL.

\lmml{071222c}
For every rainbow $\cX$, the coherent configuration $\pr_2 \WL_3(\cX)$ is sesqui\-closed.  In particular, $\WLD(\cX)\le \pr_2 \WL_3(\cX)$.
\elmm
\prf
Put $\fX=\WL_3(\cX)$; in particular, $\pr_2 \WL_3(\cX)=\pr_2\fX$. Let $\alpha\in\Omega$ and $y=(\alpha)$ a $1$-tuple. The residue $\res_y \fX$ of the $3$-ary coherent configuration~$\fX$ with respect to~$y$ is a coherent configuration on~$\Omega$, and by Lemma~\ref{100622c1} (for $m=3$),
\qtnl{071222f}
\res_y \fX\ge (\pr_2 \fX)_\alpha.
\eqtn
Let $s$ be a basis relation of the coherent configuration $\pr_2 \fX$. Then there is a class $X\in\fX$ such that $\pr_2 X=s$. Without loss of generality, we may assume that $X$ consists of the triples $(\beta,\gamma,\beta)$  with $(\beta,\gamma)\in s$.  Now if  $\alpha s\ne\varnothing$, then
$$
\{\alpha\}\times \alpha s=\pr_2 X\,\cap \,\res_y X=\res_y X.
$$
Consequently, $\{\alpha\}\times \alpha s$ is a basis relation of the coherent configuration $\res_y \fX$. It follows that $\alpha s$ is a fiber of it. Since $\alpha s$ is a homogeneity set of $(\pr_2\fX)_\alpha$, condition~(S1) follows from~\eqref{071222f}. This proves that the coherent configuration~$\pr_2 \fX$ satisfies condition~(S1).

Let $\alpha,\alpha'\in\Omega$ belong to the same fiber of  the coherent configuration~$\pr_2\fX$. Then the $1$-tuples $y=(\alpha)$ and $y'=(\alpha')$ lie in the same class of the partition~$\pr_1\fX$. By  Lemma~\ref{040622a}  for $m=3$, $k=2$, and $\varphi=\id_\fX$, the  algebraic automorphism~\eqref{110423a} % $\varphi_{yy'}: \res_{y^{}} \fX\to\res_{y'}\fX$ 
extends the identity algebraic automorphism of $\pr_2\fX$ and takes the singleton~$\{1_{\alpha^{}}\}\in\res_{y^{}}\fX$ to the singleton~$\{1_{\alpha'}\}\in\res_{y'} \fX$. Consequently, this  algebraic automorphism  has the $\alpha\alpha'$-extension. Thus the coherent configuration~$\pr_2\fX$ satisfies condition~(S2). This implies that it is sesquiclosed, and $\WLD(\cX)\le \pr_2\fX$ by the minimality of the sesquiclosure. 
\eprf

In general, not every algebraic isomorphism of two coherent configurations can be extended to an algebraic isomorphism between their sesquiclosures. An obvious reason for this is that algebraic isomorphisms do not preserve the property to be sesquiclosed. To exclude this reason, one can define sesquiclosed algebraic isomorphisms as follows~\cite{Ponomarenko2022a}: an algebraic isomorphism $\varphi:\cX\to\cX'$ is said to be  \emph{sesquiclosed} if $\varphi$ has  the $\alpha\alpha'$-extension for all  $\alpha\in\Delta\in F(\cX)$ and $\alpha'\in\Delta'\in F(\cX')$ such that $\Delta^\varphi=\Delta'$. 

We extend the concept of sesquiclosed algebraic isomorphism to rainbows by saying  that a similarity $\varphi:\cX\to\cX'$ is sesquiclosed if there exists (and hence unique) a sesquiclosed algebraic isomorphism $\psi:\WLD(\cX)\to\WLD(\cX')$ that extends~$\varphi$. In this case we also say that the rainbows $\cX$ and $\cX'$ are \emph{$\WLD$-equivalent} with respect to~$\varphi$.  The algorithm $\WLD$ \emph{identifies} a rainbow $\cX$ if every sesquiclosed similarity from~$\cX$ is induced by an isomorphism. Saying that graphs $\cG$ and $\cG'$ are $\WLD$-equivalent (respectively,  the algorithm $\WLD$ identifies $\cG$), we mean that the the rainbows $\cX(\cG)$ and $\cX(\cG')$ are $\WLD$-equivalent with respect to the standard similarity (respectively,  the algorithm $\WLD$ identifies  the rainbow $\cX(\cG)$).

\section{Sesquiclosed coherent configurations and the $3$-dim WL}\label{120223e}

%In this section we prove Theorems~\ref{050601a} and~\ref{130123b}, and Corollary~\ref{211222w}. 
Let $\cX$  be a sesquiclosed coherent configuration. We say that two point triples $x$ and $x'$ are \emph{similar} (with respect to~$\cX$) if $(x^{}_1,x^{}_2)\cdot x^{}_3\sim (x'_1,x'_2)\cdot x'_3$, i.e., if there exists the $x^{}_3x'_3$-extension $\varphi_{x^{}_3x'_3}$ of the algebraic automorphism~$\id_\cX$ and also
\qtnl{151222j}
\varphi_{x^{}_3x'_3}(r_{x^{}_3}(x^{}_1,x^{}_2))=r_{x'_3}(x'_1,x'_2).
\eqtn
It is not hard to see that the relation ``to be similar'' is an equivalence relation on the triples. 

\thrml{091222a1}
Let $\cX$ be a sesquiclosed coherent configuration on $\Omega$. Then two  triples of $\Omega^3$ are similar if and only if they belong to the same class of $\WL_3(\cX)$.  
\ethrm
\prf
Let $\cP$ be the partition of  $\Omega^3$ into the classes of similar triples, and let $\fX=\WL_3(\cX)$. We need to verify that $\cP=\fX$. To prove that $\cP\le \fX$, we  assume that $x$ and $x'$ belong to the same class $X\in \fX$.  Then the $(1)$-tuples $y=(x_3)$ and $y'=(x_3')$ are $\psi$-similar, where $\psi=\id_\fX$, and 
\qtnl{200323a}
(x_1,x_2)\in\res_{y^{}}X\qaq (x'_1,x'_2)\in\res_{y'}X.
\eqtn
Let $\psi_{yy'}:\res_{y^{}} \fX\to \res_{y'}\fX$ be the algebraic isomorphism from Lemma~\ref{040622a} for $m=3$ and $k=2$. It extends the algebraic automorphism $\psi_2=\id_{\pr_2\fX}$ and hence  the algebraic automorphism $\id_\cX$ (recall that $\pr_2\fX\ge \cX$). Furthermore,
$$
\psi_{yy'}(1_{x^{}_3})=1_{x_3'}\qaq \psi_{yy'}(\res_{y^{}}X)=\res_{y'}(\psi(X))=\res_{y'}X.
$$
Thus the restriction of $\psi_{yy'}$ to the coherent configuration $\cX_{x_3}\le \res_{y^{}} \fX$ is the $x^{}_3x'_3$-extension $\varphi_{yy'}$ of the algebraic automorphism~$\id_\cX$. By formula~\eqref{200323a}, we obtain
$$
\varphi_{yy'}(r_{y^{}}(x^{}_1,x^{}_2))=\psi_{yy'}(r_{y^{}}(x^{}_1,x^{}_2))\supseteq\psi_{yy'}(\res_{y^{}}X)=
$$
$$
\res_{y'}(\psi(X))=\res_{y'}X\ni (x'_1,x'_2).
$$
It follows that $\varphi_{yy'}(r_{y^{}}(x^{}_1,x^{}_2))=r_{y'}(x'_1,x'_2)$, i.e., the triples $x$ and $x'$ are similar. Thus, $\fX\ge \cP$.  To prove the converse inclusion, we need an auxiliary lemma.

\lmml{141222a}
Let  $\varphi:\cX\to\cX'$ be an algebraic isomorphism of sesquiclosed coherent configurations.  Assume that $(x,x')\in\Omega^3\times{\Omega'}^3$ is such that $\varphi$ has the  $yy'$-extension $\varphi_{y^{}y'}$ with $y=x_3$ and $y'=x_3'$, and also
\qtnl{060223v}
\varphi(r(y,y))=r'(y',y')\qaq \varphi_{y^{}y'}(r^{}_{y^{}}(x_1,x_2))=r'_{y'}(x'_1,x'_2).
\eqtn  
Then $(x,x')\in T^3_\varphi(\cX,\cX')$.
\elmm
\prf
Given an arbitrary  pair $(r,s)$ of basis relations of~$\cX_y$, put 
$$
\Delta(r,s)=x_1 r\cap x_2s^*.
$$ 
If this set is not empty, then the sets $x_1r$ and $x_2s^*$ are contained in the same fiber of the coherent configuration~$\cX_y$, that coincides with the right support of both~$r$ and~$s^*$; denote this fiber  by $\Gamma(r,s)$. When $\Delta(r,s)=\varnothing$, we put  $\Gamma(r,s)=\varnothing$. Similarly, we define the sets~$\Delta(r',s')$ and~$\Gamma(r',s')$ for any pair $(r',s')$ of basis relations of the coherent configuration~$\cX'_{y'}$. 

In what follows, the image of a relation $r$ with respect to $\varphi_{y^{}y'}$ is denoted by~$r'$. In these notation, the right-hand side of~\eqref{060223v} implies that if $t=r^{}_{y^{}}(x_1,x_2)$, then $t'=\varphi_{yy'}(t)=r'_{y'}(x'_1,x'_2)$, and hence
\qtnl{091233t}
|\Delta(r,s)|=c_{r,s}^{t^{}}=c_{r',s'}^{t'}=|\Delta(r',s')|.
\eqtn
Moreover, if $\Delta(r,s)\ne\varnothing$, then the  algebraic isomorphism $\varphi_{y^{}y'}$ takes the right support~$\Gamma(r,s)$ of~$r$ to the right support $\Gamma(r',s')$ of~$r'$. It follows that
\qtnl{091233t1}
\Gamma(r,s)'=\Gamma(r',s').
\eqtn

To prove that $(x,x')\in T^3_\varphi(\cX,\cX')$, we consider  the pebble game  $\sC_4(\varphi)$ with initial configuration $(x,x')$, see Subsection~\ref{060223e}. Let Spoiler choose a set $A'\subseteq\Omega'$ (the case $A'\subseteq \Omega$ is similar). Denote by $S(A')$ the set of all pairs $(r',s')$ of basis relations of~$\cX'_{y'}$, such that the set  $A'\cap\Delta(r',s')=:A'(r',s')$ is not empty. Then
$$ 
|A'|=\sum_{(r',s')\in S(A')}|A'(r',s')|.
$$
Next, for every pair $(r',s')\in S(A')$, denote by  $A(r,s)$ an arbitrary subset of $\Delta(r,s)$, the cardinality of which is equal to  the number $|A'(r',s')|\le|\Delta(r',s')|=|\Delta(r,s)|$, see~\eqref{091233t}. 

Now, Duplicator responds with the union $A\subseteq\Omega$ of all the sets $A(r,s)$ for which  $(r',s')\in  S(A')$. This is correct, because
$$
|A|=\sum_{(r',s')\in S(A')}|A(r,s)|=\sum_{(r',s')\in S(A')}|A'(r',s')|=|A'|.
$$

Next,  Spoiler picks up some $a\in A$, and Duplicator responds with arbitrary $a'\in A'(r',s')\subseteq A'$, where the pair $(r',s')\in S(A')$ is defined so that $a\in A(r,s)$. 

To complete the proof, it remains to verify that $r(x_i,a)'=r'(x'_i,a')$, $i=1,2,3$.  Note that 
$$
a\in A(r,s)\subseteq\Delta(r,s)\subseteq\Gamma(r,s)\qaq a'\in A'(r',s')\subseteq\Delta(r',s')\subseteq \Gamma(r',s'),
$$
see~\eqref{091233t1}. By condition~(S1), we have $\Gamma(r,s)=yu$ for some $u\in S$. In view of the first equality in~\eqref{060223v}, condition~(S2) implies that $\Gamma(r',s')=y'u'$. Thus, $r(x^{}_3,a)'=r(y,a)'=u'=r'(y',a')=r'(x'_3,a')$. Next, 
$$
r(x^{}_1,a)\supseteq r_y(x^{}_1,a)=r\qaq
r'(x'_1,a')\supseteq r'_{y'}(x'_1,a')=r'.
$$
Since the algebraic isomorphism $\varphi_{y^{}y'}$ coincides with $\varphi$  on $\cX$ and takes $r$ to $r'$, this shows that $r(x^{}_1,a)'\supseteq r'\ni (x'_1,a')$. Consequently, $r(x^{}_1,a)'=r'(x'_1,a')$. The same argument for $x_1$ and $r$ replaced by, respectively, $x_2$ and $s$, shows that $r(x^{}_2,a)=r'(x'_2,a')$.
\eprf

Let us complete the proof of Theorem~\ref{091222a1}, let $x$ and $x'$ be similar triples. In view of condition~(S2),  the assumption of Lemma~\ref{141222a} is satisfied for $\cX'=\cX$ and~$\varphi=\id_\cX$, and therefore $(x,x')\in T^3_\varphi(\cX,\cX')$. By Lemma~\ref{131222a}(1) for $m=3$, $\cX=\cX'$, and $\hat\varphi=\id_\fX$, this implies that $[x]=[x']$. Thus, $\cP\ge \fX$.
\eprf

\section{Proofs of the main results}\label{300323a}
In this section we first prove Theorem~\ref{050601a}, and then deduce formula~\eqref{140323b}, Theorem~\ref{130123b}, and Corollary~\ref{211222w}. Without loss of generality, we will prove all these statements for a graph $\cG$ replaced by its  coherent configuration  $\cX=\WL(\cG)$.  In what follows, we put $\fX=\WL_3(\cX)$.\medskip

{\bf Proof of Theorem~\ref{050601a}.} 
Denote by $\fY_0$  the partition of~$\Omega^3$ into the color classes of the coloring $c_0(3,\pr_2\fX)$. It is not hard to see that $\fY_0\le \fX$. Taking the $3$-ary coherent closure of both sides yields
$$
\WL_3(\pr_2\fX)=\WL_3(\fY_0)\le\WL_3(\fX)=\fX.
$$ 
Using this inclusion and applying the operator $\WL_3$ to each coherent configuration in the series $\cX\le \WLD(\cX)\le \pr_2\fX$ (see Lemma~\ref{071222c}), we obtain
\qtnl{230323e}
\fX=\WL_3(\cX)\le \WL_3(\WLD(\cX))\le\WL_3(\pr_2\fX)\le\fX,
\eqtn
whence $\WL_3(\cX)=\WL_3(\WLD(\cX))$ .

To prove the second equality, set $\cY=\WLD(\cX)$.  By Lemma~\ref{071222c}, the coherent configuration $\cY'=\pr_2\fX$ is sesquiclosed and  hence  $\cY\le \cY'$.  Assume on the contrary that $\cY<\cY'$. Then there are relations $s\in S(\cY)$ and $s'\in S(\cY')$, and a point $\alpha\in\Omega$ such that $\alpha s'$ is a proper subset of~$\alpha s$.  Take arbitrary points $\beta\in\alpha s'$ and $\beta'\in\alpha s\setminus\alpha s'$ and put 
$$
x=(\beta,\beta,\alpha) \qaq x'=(\beta',\beta',\alpha).
$$ 

The triples $x$ and $x'$ are  similar with respect to~$\cY$. Indeed, $\varphi_{\alpha\alpha}=\id_{\cY_\alpha}$ is the $\alpha\alpha$-extension of the algebraic automorphism~$\id_\cY$. Moreover,  since the coherent configuration $\cY$ is sesquiclosed, $1_{\alpha s}$ is a basis relation of its point extension~$\cY_\alpha$. Therefore, $r_\alpha(\beta,\beta)=1_{\alpha s}=r_\alpha(\beta',\beta')$ and hence 
$$
\varphi_{\alpha\alpha}(r_\alpha(\beta,\beta))=\varphi_{\alpha\alpha}(1_{\alpha s})=1_{\alpha s}=r_\alpha(\beta',\beta').
$$
By  Theorem~\ref{091222a1}, this implies that  $x$ and $x'$  lie in the same class of the $3$-ary coherent configuration~$\WL_3(\cY)$.  However, by the first part of the proof,  $\WL_3(\cY)=\WL_3(\cY')$, see formula~\eqref{230323e}. Thus, $x$ and $x'$ lie in the same class of $\WL_3(\cY')$. 

On the other hand, $r_{\cY'}(\alpha,\beta)\ne r_{\cY'}(\alpha,\beta')$. Therefore, the triples $x$ and $x'$ are not  similar with respect to~$\cY'$. Since the coherent configuration $\cY'$ is sesquiclosed, this implies by Theorem~\ref{091222a1} that $x$ and $x'$ do not lie in the same class of $\WL_3(\cY')$, a contradiction. \bull

{\bf Proof of formula~\eqref{140323b}.}  
 Denote by  $\sigma_{1,2}$ the algorithm defined in~\eqref{010123k} for the set $I=\{1,2\}$. Obviously, $\sigma_{1,2}(\cX)\le W(\cX)$.
 
 \lmml{010123u1}
 $\WLD(\cX)=\sigma_{1,2}(\cX)$.
 \elmm
 \prf
 The coherent configuration $\sigma_{1,2}(\cX)$ is sesquiclosed. Indeed, taking into account that $\sigma_1(\sigma_{1,2}(\cX))=\sigma_{1,2}(\cX)=\sigma_2(\sigma_{1,2}(\cX))$, we see that it  satisfies condition~(S1) by~\cite[Chapter~O, Theorem~4.7(b)]{Wei1976}, and satisfies  condition~(S2) by the definition of~$\sigma_1$. Thus, 
 $
 \sigma_{1,2}(\cX)\ge \WLD(\cX).
 $ 
 
 Conversely, it suffices to verify that for any two pairs $x$ and $x'$ of points of the coherent configuration~$\cY=\WLD(\cX)$, we have
 \qtnl{140323a}
 r_\cY(x)=r_\cY(x')\quad\Rightarrow\ x\simn{2}x'.
 \eqtn
 Indeed, assume that this is true. Then $\sigma_2(\cY)=\cY$. Furthermore, $\sigma_1(\cY)=\cY$ because $\cY$ is sesquiclosed (see condition~(S1). Thus,
 $$
 \sigma_{1,2}((\WLD(\cX))=\sigma_{1,2}(\cY)=\cY=\WLD(\cX)
 $$ 
 and $\WLD(\cX)=\sigma_{1,2}(\WLD(\cX))\ge \sigma_{1,2}(\cX)$, which completes the proof.

 To prove \eqref{140323a}, we assume that  $r_\cY(x)=r_\cY(x')$. Then obviously $r(x)=r(x')$. Next, let $\alpha\in\Omega$, and let $X$ be the class of $\WL_3(\cY)$ containing the triple $x\cdot\alpha$. Since $ r_\cY(x)=r_\cY(x')$ and  $\pr_2 \WL_3(\cY) =\cY$ (see above), we have $x'\in\pr_2 X$. Consequently, there exists $\alpha'\in\Omega$ such that $x'\cdot\alpha'\in X$. It follows (see Subsection~\ref{300323c}) that
 $$
 n_{\alpha^{}}(x)=n_2(X)=n_{\alpha'}(x').
 $$ 
 Thus,  $x\simn{2}x'$, as required.
 \eprf
 
Let us prove the left-hand side inclusion in formula~\eqref{140323b}.  By the second inclusion in~\eqref{230323a}, we have $\pr_2 \WL_3(\cX)=\WLD(\cX)$.  Thus  by Lemma~\ref{010123u1}, we have
$$
\pr_2 \WL_3(\cX)=\WLD(\cX)=\sigma_{1,2}(\cX)\le W(\cX),
$$
as required.

Let us prove the right-hand side inclusion in formula~\eqref{140323b}. Let $\cY'=\pr_2\WL_4(\cX)$. Since $W(\cX)$ is the minimal coherent configuration such that  $\sigma_i(\cX)=\cX$ for all~ $i$, it suffices to verify that $\sigma_i(\cY')=\cY'$ also for all $i$. By formula~\eqref{300323f}, we have $\pr_2\WL_3(\cY')=\cY'$. Consequently, by Theorem~\ref{050601a},
$$
\WLD(\cY')=\pr_2\WL_3(\cY')=\cY',
$$ 
i.e., $\cY'$ is sesquiclosed. Thus it suffices to verify that $\sigma_3(\cY')=\sigma_4(\cY')=\cY'$. Note that this is obviously true if  for any two pairs $x$ and $x'$ of points of the coherent configuration~$\cY'$, we have
\qtnl{140323a5}
r_{\cY'}(x)=r_{\cY'}(x')\quad\Rightarrow\ x\simn{3}x'\qaq x\simn{4}x'.
\eqtn

To verify this implication, assume that $r_{\cY'}(x)=r_{\cY'}(x')$. Then there are pairs $y$ and $y'$ such that in the $4$-ary coherent configuration $\fX:=\WL_4(\cX)$ we have 
\qtnl{221023a}
[y\cdot x]=[y'\cdot x'].
\eqtn
It follows that the pairs $x$ and $x'$ are $\id_\fX$-similar. By Lemma~\ref{040622a}, this implies that $(\id_\fX)_{xx'}$ induces  the $xx'$-extension of the algebraic isomorphism $\pr_2\id_\fX=\id_{\cY'}$. Thus,  $x\simn{3}x'$. 

To verify that $x\simn{4}x'$, denote by $X_1$, $X_2,\ldots,X_k$  all pairwise distinct  classes  $[\wt y\cdot x]\in \fX$ such that $x\cdot y\sim x\cdot \wt y$. Then it is easily seen that
$$
n_{y^{}}(x)=\sum_{i=1}^k|\res_x X_i|=\sum_{i=1}^k n_2(X_i).
$$
Thus to verify that $n_{y^{}}(x)=n_{y'}(x')$ (and hence $x\simn{4}x'$), it suffices to check that given a pair $\wt y$ such that $x\cdot y\sim x\cdot \wt y$, there exists a pair~$\wt y'$ such that $x'\cdot y'\sim x'\cdot \wt y'$ and the class $\wt X:=[\wt y\cdot x]$ coincides with the class $[\wt y'\cdot y']$.

We have $x\in \pr_{\{3,4\}} X\cap\pr_{\{3,4\}}\wt X$. It follows that $x'\in \pr_{\{3,4\}} X=\pr_{\{3,4\}}\wt X$. This implies that there exists a pair $\wt y'$ such that $[\wt y'\cdot y']=\wt X$.  Consequently,
$$
(\id_\fX)_{y'y}(\res_{y'}X)=\res_{y^{}} X \qaq (\id_\fX)_{\wt y\wt y'}(\res_{\wt y} \wt X)=\res_{\wt y'}\wt X.
%\res_{y'}X \xrightarrow{(\id_\fX)_{y'y}}\res_y X\qaq\res_{\wt y} \wt X \xrightarrow{(\id_\fX)_{\wt y\wt y'}}\res_{\wt y'}\wt X.
$$
Furthermore, by Lemma~\ref{100622c1}, we have  $\res_{y'}X \subseteq r_{y'}(x')$, $ \res_y X\subseteq  r_y(x)$ and $\res_{\wt y}\wt X \subseteq r_{\wt y}(x)$, $  \res_{\wt y'}\wt X\subseteq r_{\wt y'}(x')$. By Lemma~\ref{040622a}, this yields
\qtnl{221023b}
(\id_\cX)_{y'y}(r_{y'}(x'))=r_y(x)\qaq (\id_\cX)_{\wt y \wt y'}(r_{\wt y}(x) )=r_{\wt y'}(x').
%r_{y'}(x')  \xrightarrow{(\id_\cX)_{y'y}} r_y(x)\qaq r_{\wt y}(\wt x)  \xrightarrow{(\id_\cX)_{\wt y'y}}r_{\wt y'}(\wt y').
\eqtn
Finally, since  $x\cdot y\sim x\cdot \wt y$, we also have
\qtnl{221023c}
(\id_\cX)_{y\wt y}(r_y(x))=r_{\wt y}(x).
\eqtn
By \eqref{221023b} and \eqref{221023c}, the composition mapping $(\id_\cX)_{y'\wt y'}=(\id_\cX)_{y'y}\,(\id_\cX)_{y\wt y}(\id_\cX)_{\wt y \wt y'}$ takes the relation $r_{y'}(x')$ to the relation $r_{\wt y'}(x')$. Thus, $x'\cdot y'\sim x'\cdot \wt y'$, and we are done.\bull

{\bf Proof of Theorem~\ref{130123b}.} The first part of the statement follows from  Theorem~\ref{121222a} below if the similarity $\varphi$ is taken to be the standard one.

\thrml{121222a}
Rainbows $\cX$ and $\cX'$ are $\WL_3$-equivalent with respect to a similarity $\varphi:\cX\to \cX'$ if and only if  $\cX$ and~$\cX'$ are $\WLD$-equivalent  with respect to $\varphi$.
\ethrm
\prf
To prove the ``only if'' part, assume that $\cX$ and $\cX'$ are $\WL_3$-equivalent with respect to~$\varphi$, i.e., there is an algebraic isomorphism $\hat\varphi:\WL_3(\cX)\to\WL_3(\cX')$ such that  $\psi:=\hat\varphi_2$ defined  in~\eqref{150223a} for $k=2$, extends the similarity~$\varphi$. By the second equality in formula~\eqref{230323a},  the algebraic isomorphism $\psi$ takes~$\cY:=\WLD(\cX)$ to~$\cY':=\WLD(\cX')$. It remains to verify that $\psi$ is sesquiclosed. To this end, let 
$$
\alpha\in\Delta\in F(\cY)\qaq \alpha'\in\Delta'\in F(\cY')
$$ 
so that $\Delta^\psi=\Delta'$. Denote by $y$ and $y'$ the $1$-tuples $(\alpha)$ and~$(\alpha')$, respectively. By Lemma~\ref{040622a} for $m=3$ and  $k=2$, the algebraic isomorphism $\hat\varphi$ induces an algebraic isomorphism 
$$
\hat\varphi_{yy'}:\res_{y^{}}\WL_3(\cX)\to\res_{y'} \WL_3(\cX'),
$$
that extends~$\psi$.  Since the restriction of $\hat\varphi_{yy'}$ to the point extension $\cY^{}_{\alpha^{}}$, takes this extension to $\cY'_{\alpha'}$, we are done.

Let us prove the ``if'' part. Without loss of generality, we may assume that $\cX=\WLD(\cX)$ and $\cX'=\WLD(\cX')$. Let $\varphi:\cX\to\cX'$ be a sesquiclosed algebraic isomorphism. Then given  $x\in \Omega^3$, there exists $x'\in{\Omega'}\phmaa{3}$ such that condition~\eqref{060223v} is satisfied; note that this guarantees the existence of  the $yy'$-extension $\varphi_{yy'}$ of  the algebraic isomorphism~$\varphi$, where $y=(x^{}_3)$ and $y'=(x_3')$.  

We claim that the mapping 
$$
\hat\varphi:\WL_3(\cX)\to\WL_3(\cX'),\ [x] \mapsto[x'],
$$
is a correctly defined bijection. Indeed, the uniqueness of the point extension of algebraic isomorphism implies that if condition~\eqref{060223v} is satisfied for $x$ and $x'$ replaced by $\wt x\in\Omega^3$ and $\wt x'\in{\Omega'}\phmaa{3}$, respectively, then the diagram
$$
\begin{CD}
	\cX^{}_{y^{}} @> \varphi_{y^{}y'} >> \cX'_{y'} \\
	@V{(\id_\cX)_{y\wt y}} VV @VV{(\id'_{\cX'})_{y'\wt y'}}V \\
	\cX_{\wt y^{}} @>> \varphi_{\wt y^{}\wt y'}  > \cX'_{\wt y'} \\
\end{CD}\
$$
is  commutative, where $\wt y=(\wt x_3)$ and $\wt y'=(\wt x'_3)$. Here we also use the fact that $\cX$ and $\cX'$ are sesquiclosed. It follows that if $\wt x$ is $\hat\varphi$-similar to $x$, then $\wt x'$ is $\hat\varphi$-similar to $x'$. By Theorem~\ref{091222a1}, this implies that $[\wt x]=[x]$ and $[x']=[\wt x']$, which proves the claim. 

Let us complete the proof of Theorem~\ref{121222a}. By the first equality in~\eqref{230323a},  we have  $\hat\varphi_2=\varphi$. Furthermore, from Lemma~\ref{141222a}, it follows that if~$\hat\varphi([x])=[x']$ for some $x$ and~$x'$, then $(x,x')\in T^3_\varphi(\cX,\cX')$. By Lemma~\ref{131222a}(2), this shows that~$\hat\varphi$ is an algebraic isomorphism. Thus, $\cX$ and $\cX'$ are $\WL_3$-equivalent with respect to~$\varphi$.
\eprf

The second part of Theorem~\ref{130123b} follows from the first one. Indeed, the coherent configuration $\cX$ is identified by $\WL_3$ if and only if  every algebraic isomorphism with respect to which $\cX$ and $\cX'$ are $\WL_3$-equivalent is induced by an isomorphism. By Theorem~\ref{121222a}, this  is true if and only if  every algebraic isomorphism with respect to which $\cX$ and~$\cX'$ are $\WLD$-equivalent, is induced by an isomorphism, i.e., $\cX$ is identified by~$\WLD$. \bull

\section{Applications to projective planes}\label{120223r}

\subsection{Preliminaries} Let $\cX$ be a coherent configuration on~$\Omega$. Denote by $\cX\otimes\cX$ the tensor square of $\cX$, which is  a coherent configuration on~$\Omega^2$ with basis relations $r\otimes s=\{(x,y):\ (x_1,y_1)\in r,\ (x_2,y_2)\in s\}$, where $r$ and $s$ run over~$S$. 

Following~\cite{EvdP2010}, the \emph{$2$-extension} and \emph{$2$-closure} of~$\cX$ are defined to be the coherent configurations
$$
\hat\cX=\WL(\cX\otimes \cX,1_\Delta)\qaq \bar \cX=(\hat\cX_\Delta)^f,
$$ 
respectively, where $\Delta=\diag(\Omega^2)$ and $f$ is the bijection $\Delta\to\Omega$, $(\alpha,\alpha)\mapsto\alpha$. The $2$-extension accumulates information about all one point extensions of~$\cX$. More precisely, let $e$ be the equivalence relation  with classes $\Gamma=\Omega\times\{\alpha\}$, $\alpha\in\Omega$. Then $e$ is a parabolic of~$\hat\cX$, i.e., an equivalence relation on~$\Omega^2$ that is also a relation of the coherent configuration $\hat\cX$. Furthermore,
\qtnl{060123a}
(\hat\cX_\Gamma)^{f_\alpha}\ge\cX_\alpha,
\eqtn
where $f_\alpha$ is the bijection $\Gamma\to\Omega$, $(\beta,\alpha)\mapsto\beta$. 

The $2$-closure defines a closure operator which is more powerful than the operator~$\WLD$. Indeed,  from \cite[Lemma~3.5.25]{CP2019}, it follows that the coherent configuration $\bar\cX$ satisfies conditions~(S1) and~(S2), and hence is sesquiclosed. This  proves the proposition below (which strengthens the left-hand side inclusion in \cite[Theorem~1.4]{EvdP1999c}  for $m=2$).

\prpstnl{060123y}
$\WLD(\cX))\le\bar\cX$.
\eprpstn

Let $\cX'$ be a coherent configuration on $\Omega'$. An algebraic isomorphism $\hat\varphi:\wh\cX\to\wh{\cX'}$ is called a \emph{$2$-extension} of an algebraic isomorphism  $\varphi:\cX\to\cX'$ if
$$
\hat\varphi(1_\Delta)=1_{\Delta'}\qaq \hat\varphi(r\otimes s)=\varphi(r)\otimes\varphi(s)
$$
for all $r,s\in S$, where $\Delta'=\diag({\Omega'}\phmaa{2})$. Not every $\varphi$ has the $2$-extension, but if it has, then $\wh\varphi$ is uniquely determined. The following lemma is a special case of~\cite[Lemma~8.3(2)]{EvdP2000} for $m=2$.

\lmml{060123c}
If $\cX=\bar\cX$ and $\varphi$ has the $2$-extension, then $\varphi$ is sesquiclosed.
\elmm

\subsection{The scheme of a projective plane} Let $\cP$ be  a finite projective plane of order~$q$. Denote by $\Omega$ the union of the point and line sets of~$\cP$, $|\Omega|=2(q^2+q+1)$. Then the set $\Omega^2\setminus 1_\Omega$  is partitioned into three symmetric relations $s_1$, $s_2$, and $s_3$ consisting of all pairs of distinct points or lines,  of a point and a line that are incident, and  of a point and a line that are not incident, respectively. Together with $s_0=1_\Omega$ they form a symmetric scheme $\cX(\cP)$ of rank~$4$ and the valencies 
\qtnl{060123q}
n_{s_0}=1,\quad  n_{s_1}=q^2+q,\quad n_{s_2}=q+1,\quad n_{s_3}= q^2.
\eqtn
Every scheme algebraically isomorphic to $\cX(\cP)$ is of the form $\cX(\cP')$ for some projective plane $\cP'$ of the same order~$q$. Conversely, the schemes of any two projective planes of the same order are algebraically isomorphic, and, in view of~\eqref{060123q}, the corresponding algebraic isomorphism is unique.

The $2$-extension of the scheme of a projective plane of arbitrary order~$q$ was explicitly calculated in \cite[Theorem~1.4]{EvdP2010}.  It turned out that its rank does not depend on $q\ge 3$ and is equal to~$208$. Statements~(2) and~(3) of the cited theorem yield  the following result to be used in Subsection~\ref{070123a}.

\thrml{070123b}
Let $\cP$ be a projective plane and $\cX=\cX(\cP)$. Then $\bar\cX=\cX$ and every algebraic isomorphism from~$\cX$ has the $2$-extension.
\ethrm

In fact, Theorem~\ref{070123b} follows from a parametrization of the basis relations of~$\hat\cX$ by certain configurations of points and lines of~$\cP$. Using this parametrization (see~\cite[Appendix]{EvdP2010}), the basis relations contained in the parabolic~$e$ can completely be   described as in Table~\ref{fig:table}: each of them consists of only the pairs $((\beta,\alpha),(\gamma,\alpha))$, $\alpha,\beta,\gamma\in\Omega$, and the third column contains the number of basis relations satisfying the conditions in the first two columns.
\begin{table}
\begin{tabular}[h]{|c|c|c|}
\hline
$\{\alpha,\beta,\gamma\}=$   &   type   & $N$ \\
\hline
$\{x\}$         &    $r(x,x)=s_0$                         &  $1$\\ 
\hline
$\{x,y\}$      & $r(x,y)\in\{s_1,s_2,s_3\}$     &  $9$ \\
\hline
$\{x,y,z\}$  & $r(x,y)=r(y,z)=r(x,z)=s_1$       &  $2$ \\ 
\hline
$\{x,y,z\}$  & $r(x,y)=s_1$,\quad $r(y,z),r(x,z)\in\{s_2,s_3\}$       &  $2$ \\ 
\hline
\end{tabular}
\caption{$s\in S(\hat\cX)$ such that $s\subseteq e$.}
\label{fig:table}
\end{table}

\subsection{Proof of Theorem~\ref{020123i}}\label{070123a}
Let $\cG=X(\cP)$ be the incidence graph of a projective plane~$\cP$ of order~$q$, i.e., the vertex and edge sets of~$\cG$ are $\Omega$ and~$s_2$, respectively.  Thus  $\cG$ is a regular bipartite graph of degree $q+1$ and with $q^2+q+1$ vertices in each part, and also $\WL(\cG)=\cX(\cP)$. 

Let $\cP'$ be a projective plane  of order~$q$, and let $\varphi$ be the algebraic isomorphism from $\cX=\cX(\cP)$ to $\cX'=\cX(\cP')$. By Theorem~\ref{070123b}, we have $\bar \cX=\cX$, $\bar{\cX'}=\cX'$, and~$\varphi$ has the $2$-extension. By Proposition~\ref{060123y}, this shows that 
$\cX\le\bar\cX=\cX$. Hence, $\cX=\WLD(\cX)$, and similarly, $\cX'=\WLD(\cX')$. By Lemma~\ref{060123c}, this implies that $\varphi$ is sesquiclosed and hence the graphs $\cG$ and $\cG'$ are $\WLD$-equivalent.

\subsection{A one point extension}
We complete the section by calculating a one-point extension of the scheme $\cX=\cX(\cP)$. Let $\alpha\in\Omega$. We claim that
\qtnl{070123x}
\cX_\alpha=(\hat\cX_\Gamma)^{f_\alpha},
\eqtn
where $\Gamma$ and $f_\alpha$ are as in formula~\eqref{060123a}. By that formula, it suffices to verify that every basis relation~$r$ of the coherent configuration on the right-hand side of~\eqref{070123x} is a relation of the coherent configuration~$\cX_\alpha$.  Every such $r$ is of the form $(s^{}_\Gamma)^{f_\alpha}$ for some $s\in S(\hat\cX)$ contained in the parabolic~$e$. Assume first that $s$ is not one of the two relations  in the third row of   Table~\ref{fig:table}. Then
$$
r=(s^{}_\Gamma)^{f_\alpha}=(\alpha s_i\times \alpha s_j) \cap s_k:=s_{ijk}
$$
for appropriate $i,j,k\in\{1,\ldots,4\}$, and $r$ is a relation of~$\cX_\alpha$, see \cite[Lemma 3.3.5(2)]{CP2019}. 

One of the two remaining basis relations~$s$  consists of pairs $((\beta,\alpha),(\gamma,\alpha))$ for which $\alpha$, $\beta$, and $\gamma$ are collinear elements of~$\cP$ (three distinct points on the same line or three distinct lines having a common point), and if $s'$ is the another one, then 
$$
(s^{}_\Gamma)^{f_\alpha}\cup (s'_\Gamma)^{f_\alpha}=s^{}_{111}.
$$
Thus it suffices to verify that only $r=(s^{}_\Gamma)^{f_\alpha}$  is a relation of $\cX_\alpha$.  But this immediately follows from the straightforward equality 
$$
r=(s^{}_{122}\cdot s^{}_{221})\setminus 1_{\alpha s_1}.
$$  
Formula~\eqref{070123x} is completely proved. It implies that $F(\cX_\alpha)=\{\alpha s_i:\ i=0,1,2,3\}$ and 
the structure of the coherent configuration~$\cX_\alpha$ is described by the table

\def\phn#1{\phantom{$\alpha$}#1\phantom{$\alpha$}}
\begin{center}
\begin{tabular}[h]{ccccc}
\phn	& $\alpha s_0$   & $\alpha s_1$  & $\alpha s_2$ & $\alpha s_3$  \\
\end{tabular}\\
\begin{tabular}[h]{c|c|c|c|c|}
\cline{2-5}
$\alpha s_0$    & \phn{1}  &  \phn{1}  &  \phn{1}  & \phn{1} \\
\cline{2-5}
$\alpha s_1$    &    $1$  &   $3$   &   $2$ &   $2$ \\
\cline{2-5}
$\alpha s_2$    &    $1$  &   $2$   &  $2$   & $1$   \\
\cline{2-5}
$\alpha s_3$    &   $1$    &  $2$   &   $1$ &  $2$  \\
\cline{2-5}
\end{tabular}
\end{center}
\medskip
where the element in the row $\alpha s_i$ and column $\alpha s_j$ is equal to the cardinality of the set $S(\cX_\alpha)_{\alpha s_i,\alpha s_j}$.

\prpstnl{130223w}
Any one point extension of the scheme of projective plane of order~$q\ge 3$ is a coherent configuration of rank~$24$.
\eprpstn

One can also verify that the intersection numbers of the one point extension in Proposition~\ref{130223w} do not depend on~$q\ge 3$.

\providecommand{\bysame}{\leavevmode\hbox to3em{\hrulefill}\thinspace}
\providecommand{\MR}{\relax\ifhmode\unskip\space\fi MR }
\providecommand{\MRhref}[2]{%
	\href{http://www.ams.org/mathscinet-getitem?mr=#1}{#2}
}
\providecommand{\href}[2]{#2}

\section*{Acknowledgment}
The first two authors are supported by Natural Science Foundation of China (No. 12371019, No. 12161035)

\end{document}